\documentclass{amsbook}

\swapnumbers

\usepackage{amssymb,amsfonts}
\usepackage{commath}
\usepackage{verbatim}
\usepackage{mathrsfs}
\usepackage{eucal}
\usepackage[alphabetic]{amsrefs}
\usepackage{tikz}
\usetikzlibrary{matrix,arrows}

\theoremstyle{plain}   

\newtheorem{thm}{Theorem}[subsection]
\newtheorem{prop}[thm]{Proposition}
\newtheorem{cor}[thm]{Corollary}
\newtheorem{lemma}[thm]{Lemma}

\newtheorem{conj}[thm]{Conjecture}

\theoremstyle{remark}
   
\newtheorem{note}[thm]{Note}

\newtheorem{defn}[thm]{Definition}

\theoremstyle{plain}

\DeclareMathOperator{\id}{id}
\newcommand{\Sp}{{\operatorname{Sp}}}

\newcommand{\Sc}{{\operatorname{Sc}}}
\newcommand{\Sd}{{\operatorname{Sd}}}

\DeclareMathOperator{\colim}{colim}

\DeclareMathOperator{\im}{im}

\DeclareMathOperator{\El}{El}
\DeclareMathOperator{\Cat}{\mathbf{Cat}}
\newcommand{\op}{{\operatorname{op}}}

\newcommand{\Hom}{{\operatorname{Hom}}}
\newcommand{\Thetap}{\ensuremath{\Theta_{\operatorname{pcat}}}}

\newcommand{\cellset}{\ensuremath{\widehat{\Theta}}}

\newcommand{\ra}{\rightarrow}

\newcommand{\cat}[1]{{\operatorname{\mathbf{#1}}}}
\newcommand{\overcat}[2]{{(#1\downarrow #2)}}
\newcommand{\bltri}{\blacktriangleleft}

\DeclareMathOperator{\Map}{Map}
\DeclareMathOperator{\Mod}{Mod}

\DeclareMathOperator{\HMap}{H}
\DeclareMathOperator{\hocolim}{hocolim}

\newcommand{\psh}[1]{\ensuremath{\widehat{#1}}}

\newcommand{\C}{\ensuremath{\mathcal{C}}}

\newcommand{\G}{\ensuremath{\mathbb{G}}}
\newcommand{\W}{\ensuremath{\mathsf{W}}}
\newcommand{\sW}{\ensuremath{\mathcal{W}}}

\DeclareMathOperator{\heit}{ht}

\newcommand{\dfn}{\textbf}

\numberwithin{section}{chapter}

\begin{document}
\frontmatter
\title[A theory of weak $\omega$-categories]{A homotopy theory of weak $\omega$-categories}
\author[H. Gindi]{Harry Gindi}
\date{\today}
\begin{abstract} In this paper, we consider the model structure on the category of cellular sets originally conjectured by Cisinski and Joyal to give a model for the homotopy theory of weak \(\omega\)-categories.  We demonstrate first that any \(\Theta\)-localizer containing the spine inclusions \(\iota: \Sp[t] \hookrightarrow \Theta[t]\) must also contain the maps \(X\times \iota: X\times \Sp[t] \hookrightarrow X\times \Theta[t]\) for all objects \([t]\) of \(\Theta\)  and all cellular sets \(X\).  This implies in particular that a cellular set \(S\) is local with respect to the set of spine inclusions if and only if it is Cartesian-local.  However, we show that the minimal localizer containing the spine inclusions is not stable under two-point suspension, which implies that the equivalences between objects fibrant for this model structure only depend on their height-\(0\) and height-\(1\) structure.  We then try to see if adopting an approach similar to Rezk's, namely looking at all of the suspensions of the inclusion of a point into a freestanding isomorphism. We call the fibrant objects for this model structure \dfn{isostable Joyal-fibrant} cellular sets.  We understand the resulting model structure to be conjectured by a few mathematicians to give a model structure for a category of weak \(\omega\)-categories.  However, we make short work of this claim by producing an explicit example of a nontrivial contractible cofibrant strict \(\omega\)-category (with respect to the folk model structure) and showing that it is, first, not trivially fibrant, and second, proving that it is fibrant with respect to the isostable Joyal model structure.  We finally restart our approach from the beginning; bruised and battered, but also older and wiser, we construct a conjectural model structure which appears to have all of our desired properties.  
\end{abstract}

\maketitle

\tableofcontents

\chapter*{Introduction}
Based on the earlier work of Boardman and Vogt, Joyal introduced the theory of quasicategories as a way of dealing with constructions in homotopy theory that were modeled by homotopy-coherent diagrams.  Jacob Lurie then extended much of ordinary category theory to these quasicategories, which he and others have envisioned as a theory of \((\infty,1)\)-categories. However, as the theory of \((\infty,1)\)-categories became more and more well-understood, it became clear that it was perhaps less interesting than had been previously suspected, One major problem with the theory is that much of the richness we would expect in a theory of higher categories was simply not present. Limits and colimits are defined using cocones and cones.  Arrows that are equivalences are parameterized by maps out of the nerve of a freestanding isomorphism.  We have no notions of laxness, and we even have the property that every weighted limit or colimit can be obtained as the conical limit of an easily modified diagram.   

Moreover, since the notion of an \((\infty,2)\)-category is not a priori defined, we can't perform (homotopy coherent) 2-categorical operations on the category of \((\infty,1)\)-categories, which forces us to take the ``local'' view  of important things like Kan extensions, adjunctions, representability, and things that are generally understood using strict 2-categories.  Interesting and important higher-categorical algebraic structures like monads, operads, monoidal products, braiding, symmetry, etc. are shoehorned into what amounts to a ``thickened'' theory of \(1\)-categories, and these kinds of maneuvers make much of the theory unwieldy. 

We give an example of this: To define the Grothendieck correspondence for diagrams of spaces indexed by an \((\infty,1)\)-category, Lurie first obtains a rather large family of simplicial Quillen adjunctions \[\operatorname{St}_\phi:(\psh{\Delta}\downarrow X)_{\operatorname{Dex}} \rightleftarrows  (\psh{\Delta}^{\mathcal{X}})_{\operatorname{Proj}}\] indexed by triples \((X, \phi, \mathcal{X})\) consisting of a simplicial set \(X\), a simplicially enriched category \(\mathcal{X}\), and a map \(\phi:\mathfrak{C}[X]\to \mathcal{X}\).  Lurie shows in \cite{htt} that when this map is a weak equivalence of simplicially enriched categories, the image of the simplicially enriched subjecategory of \((\psh{\Delta}^{\mathcal{X}})_{\operatorname{Proj}}\) spanned by the cofibrant-fibrant objects under the right adjoint \((\operatorname{St}_\phi,\operatorname{Un}_\phi)\) is weakly equivalent as a simplicially enriched category to the simplicially-enriched subcategory of  \(\psh{\Delta}\downarrow X)_{\operatorname{Dex}} \) spanned by its fibrant objects.  Then by a theorem regarding simplicial model categories, this implies that the Quillen adjunction is a Quillen equivalence.  

However, from the point of view of a theory of \((\infty,2)\)-categories, we never really need to leave to an external theory in order to formulate the \((\infty,1)\)-categorical Grothendieck construction.  In fact, we have a number of different ways that we can do it.  One way to do it is to find a fibrant object that models the  \((\infty,1)\)-category of spaces and then either define the grothendieck construction by means of an oplax colimit or by means of the pullback of the generalized universal fibration, which can easily be defined using an (op)lax version of the join construction. However, this leaves us with the problem of how to deal with higher and higher dimensional versions of \((\infty, n)\)-categories to be able to easily perform \(n+1\)-dimensional operations on \((\infty,n)\)-categories.  We should think about the grothendieck construction as encoding \(n+1\)-dimensional data as \(n\)-dimensional data.  

Then to be able to formalize everything in a single theory, we're forced to consider the case where \(n\) is infinite.  In this paper, we will consider a number of candidates for a homotopy theory that encodes these ideas.  Ultirmately, we will show that the model structure proposed by Cisinski and Joyal doesn't quite pass muster, since it does not satisfy the expected stability property (namely that the two-point suspension of a Cisinski-Joyal weak equivalence need not be a weak equivalence).  However, we can use some of the tools developed in the first chapter to obtain a ``stabilized'' form of this theory.  Moreover, we show that both the original form and its stabilized counterpart are cartesian-closed model structures.  That is, it is a theorem of Cisinski that the minimal Cisniski model structure on a presheaf category is generated by choosing the cylinder functor to be the cartesian product with the subobject classifier.  This model structure is obviously cartesian-closed, by construction.  We then include the spine inclusion maps as weak equivalences, which generate the Cisinski-Joyal model structure and show that these also generate a cartesian model structure.   Finally, the stabilization of this model structure under the categorical two-point suspension functor is proven to be be cartesian closed.

However, both of these efforts are proven to fail, since we find that two strict \(\omega\)-categories that are equivalent in the sense of \cite{lmw} need not have weakly equivalent nerves in either of these model structures.  From these two failures, we see what the ultimate model structure is forced to be and state its definition in terms of generating maps, but our resulting set of generators is extremely unparsimonious.  We conjecture \eqref{genconjecture} that a smaller set of generators will suffice. 

We recently became aware of independent work reaching similar conclusions about the conjecture of Cisinski and Joyal through somewhat different methods by Dimitri Ara \cite{ara-new}.  While we provide a sequence of counterexamples in each higher dimension, he proves very generally that the suspension of any map belonging to a general class of weak equivalences in the \(\infty,n\)-case of the conjecture fails to a weak equivalence, and to rectify this, echoing Rezk, he proposes a small family of generating weak equivalences in order to rectify the problem.  

The author would like to give special thanks to Denis-Charles Cisinski for his invaluable guidance and mentorship, as well as for his help in formulating parts of \eqref{sec:cishelp}, to David Oury for his lively correspondence and his continued interest in our work, and finally, to Richard Steiner, whose tireless work in the area of pasting theory has given us a better understanding of the lax tensor product than we ever even thought was reasonable.

\mainmatter
\chapter{Background}
\section{Iterated wreath products and the category $\Theta$}
\label{sec:wreath}

We will make use of two equivalent definitions of the category \(\Theta\) of cell objects: The first definition, covered in this section is due to Berger in \cite{berger-iterated-wreath}, where he defines it to be the filtered union of wreath powers of the simplex category \(\Delta\) along the inclusion maps \(\Delta^{\wr n}\cong \Delta^{\wr n}\wr \ast \hookrightarrow \Delta^{\wr n+1}\).  This section will give a quick review of this theory.  Some of the exposition in this section is based on \cite{cisinski-decalage}, and the author makes no claims of originality in this section. 

\subsection{The category $\Gamma$}

Recall that Segal's category, \(\Gamma\), is defined as follows: The objects are the (possibly empty) sets \(\Gamma_m=\{x\in \mathbf{N}\colon 1\leq x\leq m\}\) for each \(m\in \mathbf{N}\), and morphisms \(\Gamma_m\to \Gamma_n\) are functions \(f\colon \Gamma_m\to 2^{\Gamma_m}\) such that \(f(a)\cap f(b) = \emptyset\) if \(a<b\)).  Given morphisms, \(f\colon \Gamma_m\to \Gamma_n\) and \(g\colon \Gamma_p\to \Gamma_m\), we define the composite \(f\circ g  \colon \Gamma_p\to \Gamma_n\) by letting \((f\circ g)(s)=\bigcup_{t\in f(s)}g(t)\).  It is left as an exercise to the reader to show that this law of composition is indeed associative.  

\begin{prop}The category \(\Gamma\) is equivalent to the category \(\overcat{0}{\cat{Fin}}^\op\), where \(\cat{Fin}\) is defined to be the category of finite sets, and where \(0\) denotes the set with a single element.
\end{prop}
\begin{proof}
It is immediate that \(\overcat{0}{\cat{Fin}}\) is equivalent to the subcategory of \(\cat{Set}\) whose objects are the sets \(n=\{x\in \mathbf{N}\colon 0\leq x\leq n\}\) and whose morphisms are those functions \(f\colon n\to m\) such that \(f(0)=0\).  

Given such a function, we define a map \(\Gamma_f:\Gamma_m\to \Gamma_n\) by the rule \(x\mapsto f^{-1}(x)\).  Conversely, given a map \(f\colon \Gamma_m\to \Gamma_n\), we define the subset \(f(0)\subseteq \Gamma_n\) to be the complement \[f(0)=\Gamma_n - \bigcup_{s\in \Gamma_m} f(s).\]  

Then we define the function \([f]:n\to m\) by the rule \([f](s)=j\) where \(j\) is the unique number \(0\leq j\leq m\) such that \(s\in f(j)\).  It is clear that\([\Gamma_f]=f\) and \(\Gamma_{[g]}=g\), so this determines an anti-equivalence of categories.
\end{proof}

\subsection{The categorical wreath product}

Let \(A\) be a category, and let \(F\colon B\to \Gamma\) be an object of \(\overcat{\cat{Cat}}{\Gamma}\).  

Then we define the \dfn{wreath product} \(B\wr A\) as follows: The objects of \(B\wr A\) are pairs \((b, \{a_i\}_{i\in F(b)})\) comprising an object \(b\) of \(B\) and a family of objects of \(A\) indexed by the elements of \(F(b)\).  

A morphism \((b, \{a_i\}_{i\in F(b)})\to (b', \{a'_i\}_{i\in F(b')})\) is given by the data of a pair \((f, \{\eta_{ij}\})\) comprising
\begin{enumerate}
\item [(i)] a morphism \(f \colon b\to b'\) of \(B\), and
\item [(ii)] a morphism \(\eta_{ij}:c_i\to d_j\) for each pair \(i,j\) such that \(i\in F(b)\) and \(j\in F(f)(i)\)
\end{enumerate}

It is left as an easy exercise to show that the composition of two such maps obtained in the obvious way is indeed associative.

If \(G\colon (B',F')\to (B,F)\) is a functor over \(\Gamma\), and \(\Psi\colon A'\to A\) is any functor, we obtain a functor \(G\wr F\colon B' \wr A'\to B\wr A\) defined on objects by the rule \[(b',(a'_i)_{i\in F'(b')})\mapsto (G(b),(\Psi(a'_i))_{i\in F(G(b'))}\] (which makes sense since \(F\circ G=F'\)) and is defined on morphisms in the obvious way. 

\begin{defn} We say that a category \(C\) is \dfn{semi-additive} if it admits finite products and a null object.  A morphism between such categories is a finite-limit preserving functor
\end{defn}

\begin{prop}The category \(\Gamma\wr A\) is the free semi-additive category on A, and there exists a canonical functor \(\alpha:\Gamma\wr\Gamma \to \Gamma\) sending an n-tuple of objects \(\Gamma_n(\Gamma_{i_1},\dots \Gamma_{i_n})\) to their sum, \(\Gamma_{i_1+\dots+i_n}\).  
\end{prop}
\begin{proof} 
Since \(\Gamma_0\) is a null object for \(\Gamma\), we see that \(\Gamma\wr A\) has null object equal to \(\Gamma_0()\).  We also have for any two objects \(X=\Gamma_i(a_1,\dots a_i)\) and \(Y=\Gamma_j(b_1,\dots b_j)\) an object \(\Gamma_{i+j}(a_1,\dots,a_i,b_1,\dots b_j)\), which, when equipped with the two projections \(\Gamma_{i+j}(a_1,\dots,a_i,b_1,\dots b_j)\to \Gamma_i(a_1,\dots,a_n\) and \(\Gamma_{i+j}(a_1,\dots,a_i,b_1,\dots b_j)\to \Gamma_j(b_1,\dots,b_j\) is easily checked to be a cartesian product of \(X\) and \(Y\).  It also admits a unique embedding \(A\hookrightarrow \Gamma\wr A\) defined on objects by the rule \(a\mapsto \Gamma_1(a)\), so every nonzero object of \(\Gamma\wr A\) is uniquely a product of objects of the form \(\Gamma_1(a)\). 

Since \(\Gamma\cong \Gamma\wr \ast\), and \(\Gamma\wr \ast\) is semi-additive, there exists a unique finitely continuous functor \(\Gamma\wr \Gamma\to\Gamma\wr \ast \cong \Gamma\), which sends an object 
\(\Gamma_n(\Gamma_{i_1},\dots \Gamma_{i_n}) = \prod_{j=1}^n\Gamma_1(\Gamma_{i_j})\) 
to the object \(\prod_{j=1}\Gamma_{i_j}=\Gamma_{\sum_{j=1}^n i_j}\).  
\end{proof}

\begin{prop} The wreath product is a monoidal product for the category \(\overcat{\cat{Cat}}{\Gamma}\) with monoidal unit given by the functor classifying \(\Gamma_1\), \(e_{\Gamma_1}:*\to \Gamma\).
\end{prop}
\begin{proof}
We see that \((*,e_{\Gamma_1})\wr (B,\lambda_B) \to \Gamma\wr\Gamma\) sends the objects \(\ast(b)\) to the objects \(\Gamma_1(\lambda_B(b))\), which maps under \(\alpha\) to the object \(\Gamma_{\lambda_B(b)}\) of \(\Gamma\).  Similarly 
\((B,\lambda_B)\wr (*,e_{\Gamma_1})\) sends the objects \(b(\ast,\ast,\dots,\ast)\) to the objects \(\lambda_B(b)(\ast,\ast,\dots,\ast),\) which map under \(\alpha\) to the objects \(\Gamma_{\sum_{i=1}^{\lambda_B(b)} 1}\), which are precisely the objects \(\Gamma_{\lambda_B(b)}.\)

Let \((A,f_A)\), \((B,f_B)\), \((C,f_C)\) be categories over \(\Gamma\) (we will suppress the functors \(f_X\) unless otherwise noted) to.  To prove the associativity of \(\wr\), we see that there is an isomorphism of categories, natural in \(A,B,C\), \(\alpha_{ABC}:(A\wr B)\wr C\to A \wr (B\wr C) \) where the object 
\((a, \{b_i\}_{i\in f_A(a)})(\{c_j: j\in F_b(f_A(a))\})\) is sent to the object 
\((a,\{b_i, \{c_j\}_{j\in f_B(b_i)}\}_{i\in f_A(a)}).\)

The definition on morphisms can be extracted from the definition on objects by reindexing, and we leave an explicit description of this reindexing as an exercise.  Naturality in \(A,B,C\) follows from the functoriality of the wreath product. From this, we obtain a commutative square in \(\cat{Cat}\) 

\begin{equation*}
\begin{tikzpicture}
\matrix (b) [matrix of math nodes, row sep=3em,
column sep=3em, text height=1.5ex, text depth=0.25ex]
{ (A\wr B)\wr C & A\wr(B\wr C) \\
   (\Gamma\wr \Gamma)\wr \Gamma &  \Gamma\wr(\Gamma\wr\Gamma) \\};
\path[->, font=\scriptsize]
(b-1-1) edge (b-1-2)
        edge (b-2-1)
(b-2-1) edge (b-2-2)
(b-1-2) edge (b-2-2);
\end{tikzpicture},
\end{equation*}
from which it follows that \(\alpha_{ABC}\) is indeed a morphism over \(\Gamma\) for all triples \(A,B,C\) if and only if the isomorphism \(\alpha_{\Gamma\Gamma\Gamma}\) is a morphism over \(\Gamma\).  However, it is easy to see that this holds, ultimately, by the generalized associativity of iterated addition in \(\mathbf{N}\).  
\end{proof}

\subsection{Infinite wreath products}\label{infwreath}

Iterating the wreath product construction on a category \(F_B:B\to \Gamma\) over \(\Gamma\), we obtain by recursion a definition of the \(n^\mathrm{th}\) wreath power \(B^{\wr n+1}=B \wr B^{\wr n}\).  Suppose further that \(B\) is equipped with a functor \(e_b\colon \ast \to B\) classifying an object \(b\) of \(B\). Then by iterating the wreath product construction on \(e_b:\ast\to B\), we construct the following data by recursion:
\begin{enumerate}
\item[(i)] Let \(B_0= \ast\), and let \(\iota_0=e_b\).
\item[(ii)] Let \(B_{i+1}=B\wr B_i\), and let \(\iota_{i+1}=\id_{B}\wr \iota_{i}:B\wr B_{i}\to B\wr B_{i+1}\).
\end{enumerate}
Since \(B\wr \ast\) is canonically isomorphic to \(B\), we obtain a diagram, \(T_{B,b,F_B}:\mathbf{N}\to \cat{Cat}\):
\begin{equation*}
\begin{tikzpicture}
\matrix (a) [matrix of math nodes, row sep=3em,
column sep=3em, text height=1.5ex, text depth=0.25ex]
{ B_0 & B_1 & B_2 & \dots & B_n & B_{n+1} & \dots \\};
\path[right hook->, font=\scriptsize]
(a-1-1) edge node[auto]{\(\scriptstyle \iota_0\)} (a-1-2)
(a-1-2) edge node[auto]{\(\scriptstyle \iota_1\)} (a-1-3)
(a-1-3) edge node[auto]{\(\scriptstyle \iota_2\)} (a-1-4)
(a-1-4) edge node[auto]{\(\scriptstyle \iota_{n-1}\)} (a-1-5)
(a-1-5) edge node[auto]{\(\scriptstyle \iota_n\)} (a-1-6)
(a-1-6) edge node[auto]{\(\scriptstyle \iota_{n+1}\)} (a-1-7);
\end{tikzpicture}
\end{equation*}
 
We then define \(C(B,b,F_B)=C(B,b)=\varinjlim T_{B,b,F_B}\) as the colimit of this system.  

\subsection{The simplex category $\Delta$}

Recall that the simplex category, \(\Delta\), is defined to be the skeleton of the full subcategory of \(\cat{Cat}\) spanned by the finite nonempty
linearly-ordered sets (regarded as categories).  The objects of \(\Delta\) are isomorphism classes of linearly ordered sets, where \([n]\) denotes the
class of the linearly-ordered set \(\{0<\dots<n\}\).  In fact, we may identify the skeleton with the full subcategory spanned by such sets.  In the sequel,
we will make this identification without timidity, justified by the fact that there is at most one isomorphism between any two linearly-ordered sets.

Following Rezk in \cite{rezk-theta-n-spaces}, we call a map \(f \colon [n]\to [m]\) in \(\Delta\) \dfn{sequential} if 
\[
f(i-1)+1 \geq f(i),\qquad \text{\(1\leq i \leq n\)}.
\]

We say that an object \(\gamma:[1]\to [n]\) of \(\overcat{[1]}{\Delta}\) is an interval if \(\gamma(0)=0\) and \(\gamma(1)=n\), and we say that an interval is strict if the map \(\gamma\) is also injective.  Let \(\mathcal{D}^1\) denote the full subcategory of \(\overcat{[1]}{\Delta}\) spanned by the strict intervals.  We denote a strict interval whose underlying simplex is \([n]\) by \(|n|\), and we denote the image of the inclusion \(\gamma\) by \(\partial |n|\).  

We have a functor \(q:\mathcal{D}^1\to \overcat{0}{\cat{Fin}}\) defined on objects by the formula \[|n+1|\mapsto |n+1|/\partial|n+1|=n\] and defined on morphisms by the universal property of quotients.  This gives us a functor \(q^\op:  (\mathcal{D}^1)^\op \to \overcat{0}{\cat{Fin}}^\op\cong \Gamma\).  However, the category \(\mathcal{D}^1\) is isomorphic to \(\Delta^{op}\) by the functor \([n]\mapsto \Delta([n],[1])=|n+1|\), the confirmation of which we leave to the reader.  

\label{segfun}We write \(F_\Delta\) for the induced functor \(\Delta\to \Gamma\).  We see that clearly, \(F_\Delta([m])=\Gamma_m\), and given \(f:[n]\to [m]\), we compute \(F_\Delta(f)(i)\) for \(i\in \Gamma_n\).  First, we obtain a morphism \[f^\ast:|m+1|=\Delta([m],[1])\to \Delta([n],[1])=|n+1|,\] which descends to a morphism \(q(f^\ast):n=q(|n+1|)\to q(|m+1|)=m\), so \(F_\Delta(f)(i)=\Gamma_{q(f^\ast)}(i)=(q(f^\ast))^{-1}(i)\), but \(q(f^{\ast})(j) = i\) if and only if \(f^\ast(j)=i\).  Let \(c_j:[m]\to [1]\) be the unique morphism such that \(j=\inf(c_j^{-1}(\{1\})\).  Then \(f^\ast(j)=i\) holds if and only if \(i=\inf((c_j\circ f)^{-1}(1))\) if and only if \(f(i-1)<j\leq f(i)\).  Then \(F_\Delta(f)(i)=\{j: f(i-1)<j\leq f(i)\}\).  This gives us an explicit description of the functor \(F_\Delta:\Delta\to \Gamma\).  Combining this with the definition of the wreath product, we obtain:

\subsection{The category $\Delta \wr \C$}

Let \(\C\) be a category.  Then applying the wreath product construction with respect to the functor \(p\colon\Delta\to \Gamma\), we may describe the category \(\Delta \wr \C\) explicitly as follows: An object of \(\Delta \wr \C\) is given by the data of a pair \((n,(c_1,\dots c_n))\), written \([n](c_1,\dots,c_n\), 
where \(n\in \mathbf{N}\) and \((c_1,\dots,c_n) \in Ob(\C^{\times n})\).  

A morphism \([n](c_1,\dots,c_n)\to
[m](d_1,\dots,d_m)\) is given by the data of a pair \((f, \{\eta_{ij}\})\) comprising
\begin{enumerate}
\item [(i)] a morphism \(f \colon [n]\ra [m]\) of \(\Delta\), and
\item [(ii)] a morphism \(\eta_{ij}:c_i\to d_j\) for each pair \(i,j\) such that \(f(i-1)<j\leq f(i)\)
\end{enumerate}

In general, for any category \(\C\), we will call the category \(\Delta \wr \C\) the \dfn{\(\Delta\)-suspension} of \(\C\).  
We will define \(\Theta\) to be \(C(\Delta,[0])\). 

\section{Strong generators and completions}

We will find it extremely useful to sharpen Cisinski's theory \cite{cisinski-book}*{1.4}with respect to how localizers are generated with respect to simplicial completions and how to deal with regularity \cite{cisinski-book}*{3.4} with these generators.

\subsection{Simplicial generators for localizers}

Let \(\C\) be a small category.  We let \(\mathsf{W}_\infty\) denote the \(\C\times \Delta\)-localizer generated by the maps \(X\times \Delta_n\to X\times \Delta_0\) for every presheaf \(X\) on \(\C\) and every \(n\geq 0\).  
\begin{prop}[\cite{cisinski-book}*{Corollary 2.3.7}] The \(\C\times \Delta\)-localizer \(\mathsf{W}_\infty\) is accessible.
\end{prop}   
\begin{proof}See the proof in \cite{cisinski-book}.
\end{proof}
\begin{defn}
We say that a class of maps \(S\) in \(\psh{\C}\) is a weak class of irregular generators for a \(\C\)-localizer \(\W\) if \(\W(S)=\W\).

We say that a class of maps \(S\) in \(\psh{\C\times \Delta}\) is a class of \dfn{simplicial irregular generators} for a localizer \(\W\) if the \(\C\times \Delta\)-localizer \[\mathsf{W}(S\times \Delta_0 \cup \mathsf{W}_\infty)\] is exactly the simplicial completion of \(\W\).  

We say that a class of maps \(S\)  in \(\psh{\C}\) is a class of \dfn{strong irregular generators} for \(\W\)  if the class \(S\times \Delta_0\) of maps of the form \(f\times \Delta_0\) where \(f\in S\) is a class of simplicial irregular generators for a \(\C\)-localizer \(\W\).
\end{defn}
\begin{prop}If \(S\) is a class of strong irregular generators for a \(\C\)-localizer \(\mathsf{W}\), then \(\W=\W(S)\).
\end{prop}
\begin{proof}This follows immediately from \cite{cisinski-book}*{Proposition 2.3.30}.  
\end{proof}
There is a useful and na\"ive way to strengthen classes of weak irregular generators to classes of strong irregular generators:
\begin{prop} If \(S\) is a class of weak irregular generators for a localizer \(\W\), then \(S\cup \operatorname{cart}(\{\ell:L\to e\})\) is a strong class of generators, where \(L\) is the subobject classifier of \(\C\), and \(\operatorname{cart}(\{\ell\})\) is the class of all maps \(X\times \ell:X\times L \to X\) where \(X\) is a presheaf on \(\C\).  
\end{prop}
\begin{proof}This again follows immediately from \cite{cisinski-book}*{Proposition 2.3.30}.
\end{proof}
\subsection{Strong regular generators for regular localizers}
To utilize this notion of strong generation for a localizer in the context of regular localizers, the following important proposition will be extremely useful:
\begin{prop}If \(S\) is a class of strong simplicial irregular generators for a \(\C\)-localizer \(\mathsf{W}\), then the simplicial completion of the regular completion \(\mathsf{R}(\W)\) of \(\W\) is the \(\C\times \Delta\)-localizer \(\W(S \cup \operatorname{R}(\W_\infty))\), where \(\mathsf{R}(\W_\infty)\) is the regular completion of \(\W_\infty\), which is precisely the class of objectwise weak homotopy equivalences of simplicial presheaves.
\end{prop}
\begin{proof}This follows easily from \cite{cisinski-book}*{Corollary 3.4.47}.
\end{proof}
Based on this proposition, we can give a slightly weaker notion of strong generation:
\begin{defn} 
We say that a class of maps \(S\) in \(\psh{\C}\) is a class of \dfn{weak regular generators} for a regular localizer \(\W\) if \(S\) is a class of weak irregular generators for some \(\C\)-localizer \(\W^\prime\) whose regular completion \(\mathsf{R}(\W^\prime)\) is exactly \(\W\). 

We say that a class of maps \(S\) in \(\psh{\C\times \Delta}\) is a class of \dfn{simplicial regular generators} for a regular localizer \(\W\) if \(S\) is a class of simplicial irregular generators for some \(\C\)-localizer \(\W^\prime\) whose regular completion \(\mathsf{R}(\W^\prime)\) is exactly \(\W\). 

We say that a class of maps \(S\)  in \(\psh{\C}\) is a class of \dfn{strong regular generators} for \(\W\)  if the class \(S\times \Delta_0\) of maps of the form \(f\times \Delta_0\) where \(f\in S\) is a class of simplicial regular generators for a \(\C\)-localizer \(\W\).

Unless otherwise noted, when \(\W\) is a regular localizer, a class of \dfn{strong generators} for \(W\) will always mean a class of strong \emph{regular} generators.   
\end{defn}
Then we easily obtain the following useful corollary:
\begin{cor} If \(S\) is a small set of strong generators for a regular localizer \(\W\), then the simplicial completion of \(\W\) is the class of weak equivalences of the left Bousfield localization of \(\psh{\C\times \Delta}_{\operatorname{inj}}\) at \(S\).   
\end{cor}
\begin{prop}\label{boostgenerators} If \(S\) is a class of weak regular generators for a localizer \(\W\), then \(S\cup \operatorname{cart}(\{\ell:L\to e\})\) is a strong class of regular generators, where \(L\) is the subobject classifier of \(\C\), and \(\operatorname{cart}(\{\ell\})\) is the class of all maps \(X\times \ell:X\times L \to X\) where \(X\) is a presheaf on \(\C\).  
\end{prop}
\begin{proof}This again follows immediately from from \cite{cisinski-book}*{Corollary 3.4.47}.
\end{proof}

\section{A model structure on $\widehat{\Delta \wr \C}$}\label{sec:cishelp}
We thank Denis-Charles Cisinski for his invaluable help with the formulation of this section. We will give a model structure whose fibrant objects are models for categories weakly enriched in the homotopy theory of \(\W\)-fibrant presheaves of sets on \(\C\) whenever \((\C,\W)\) a pair comprising a small category \(\C\) together with a fixed accessible cartesian regular \(\C\)-localizer, \(\W\).  For now, we fix the small category \(\C\).

\subsection{The intertwining functor $V_{\C}$}
For any category \(A\), we let \(Y_A:A\hookrightarrow \psh{A}\) denote the Yoneda embedding. Then we have an apparent pair of functors \[Y_{\Delta \wr \C}:\Delta \wr \C \hookrightarrow \psh{\C},\] the Yoneda embedding of \(\Delta \wr \C\), and by the functoriality of the wreath product, the \(\Delta\)-suspended Yoneda embedding of \(\C\), \[L=\id_\Delta\wr Y_\C:\Delta \wr \C \hookrightarrow \Delta \wr \psh{\C}\].  

We define the \dfn{\(\C\)-intertwiner} \(V_\C:\Delta \wr \psh{\C} \to \psh{\Delta \wr \C}\) to be the left Kan extension \(L_!(Y_{\Delta \wr \C})\) of \(Y_\C\) along \(L\).  Unless there is a risk of confusion, we will typically suppress the subscript \(\C\).

\subsection{Mapping objects}

For any \(\psh{\C}\)-enriched simplicial set \(X\), equipped with a pair of vertices \((x_0,x_1)\) of \(X,\) we will construct a mapping object \(X(x_0,x_1)\) of \(\psh{\C}\).

The following lemma is due to Rezk in \cite{rezk-theta-n-spaces}:
\begin{lemma}\label{leftadjointness}
Given any two families \(A_1,\dots, A_m\) and \(B_1,\dots, B_n\) of presheaves on \(\C\), the functor \(P:\psh{\C}\to \psh{\Delta\wr \C}\) defined by the formula \[X\mapsto V[n+1+m](A_1,\dots,A_m,X,B_1,\dots,B_n)\] is a parametric left adjoint, that is to say, the natural factorization \[P_0:\psh{\C}\to \overcat{P(\emptyset)}{\widehat{\Delta\wr \C}}\] of \(P\) through the forgetful functor \[U_0: \overcat{P(\emptyset)}{\psh{\Delta\wr \C}}\to \psh{\Delta\wr \C}\] admits a right adjoint. Further, we have that \[P(\emptyset)=V[m](A_1,\dots,A_m)\coprod V[n](B_1,\dots,B_m)\].   
\end{lemma}
\begin{proof}
Since we are taking the left Kan extension of the Yoneda embedding along \(L=\id_\Delta\wr Y_\C\), if we let  \(h_Z\), for any object \(Z\) of \(\Delta \wr \psh{\C}\), be the functor \(A\mapsto \Hom_{\Delta\wr\C}(A,Z)\) representing \(Z\),  we  obtain a simple formula for \(VZ\) as \(L^\ast(h_X)\) because the conical formula for the pointwise left Kan extension degenerates on the Yoneda embedding.

To see why this is true, notice that in the conical formula for the left Kan extension, we have that \[V(Z)=\varinjlim(\overcat{L}{Z}\to \Delta\wr\C\to \psh{\Delta\wr\C}),\] where \(\overcat{L}{Z}\) is the pullback \(\Delta\wr \C \to \Delta\wr\psh{\C} \leftarrow \overcat{\Delta\wr\psh{\C}}{Z}\).  However, by inspection, the category \(\overcat{L}{Z}\) is precisely the category of elements of the \(\Delta\wr\C\)-presheaf \(L^\ast (h_Z)\), so composing this diagram with the Yoneda embedding and taking a colimit is precisely the colimit of the category of elements of the presheaf \(L^\ast(h_Z)\), which just so happens to be \(L^\ast (h_Z)\) by Yoneda's lemma.  

Let \(a=[q](c_1,\dots,c_q\) be an object of \(\Delta\wr\C\).  Following Rezk in \cite{rezk-theta-n-spaces}, we see that the set of maps \(a\to L(X)\) belongs, can be divided into partitions corresponding to the partitions of \(\Hom_\Delta([q],[m+1+n])\), parameterized by the elements \(p\in \Hom_\Delta([q],[1])=\set{p}{0\leq p\leq q+1}\) as follows: 
\[G(p)=
\begin{cases}
\set{\delta}{\delta(0)\geq m+1}\qquad \text{if \(p=0\)}\\
\set{\delta}{\delta(p-1)\leq m, \delta(p)\geq m+1} \qquad \text{if \(1\leq p\leq q\)}\\
\set{\delta}{\delta(q+1)\leq m}\qquad \text{if \(p=q+1\)}
\end{cases},
\]
which decomposes the set \(\Hom_{\cellset}(a,L(X))\) into the factors \((S_0,\dots,S_{q+1})\), where the factor \(S_0\) is 
\[
\coprod_{\delta\in G(0)}\,
\prod_{i=1}^q\,\prod_{j=\delta(i-1)+1}^{\delta(i)} B_{j-(m+1)}(c_i) \approx
V[n](B_1,\dots,B_n)(\theta),
\]
the factor \(S_{q+1}\) is 
\[
\coprod_{\delta\in G(q+1)}\,
\prod_{i=1}^q\, \prod_{j=\delta(i-1)+1}^{\delta(i)} A_j(c_i) \approx
  V[m](A_1,\dots, A_m)(\theta),
\]
and the factor \(S_p\) for \(1\leq p \leq q\) is 
\[
\coprod_{\delta\in G(p)}
\left(\prod_{i=1}^{p} \,\prod_{j=\delta(i-1)+1}^{\min(\delta(i),m)} A_j(c_i)
\right) \times X(c_p) \times \left(\prod_{i=p}^q \,
  \prod_{j=\max(\delta(i-1),m)+2}^{\delta(i)} B_{j-(m+1)}(c_i)\right).
\]

It follows by inspection that the functor \(P_0\) preserves colimits and that \[P(\emptyset)=V[m](A_1,\dots,A_m)\coprod V[n](B_1,\dots,B_n).\]
\end{proof}
Since \(V[1](\emptyset)=\ast\coprod \ast\), the preceding lemma in the case where \(m=n=0\) gives us our desired right adjoint \(R:\overcat{V[1](\emptyset)}{\psh{\Delta\wr\C}}\to \psh{\C}\).  Given a \(\psh{C}\)-enriched simplicial set \(X\) together with a pair of vertices \((x_0,x_1)\) of \(X,\) we can take these data together to give a map \((x_0,x_1):V[1](\emptyset)\to X\), which give an object \(X,(x_0,x_1)\) of \(\overcat{V[1](\emptyset)}{\psh{\Delta\wr\C}}\).  Then we define \(X(x_0,x_1)=R(X,(x_0,x_1))\).  By functoriality, for any map \(f:X\to Y\) in \(\psh{\Delta\wr \C}\) and any pair of vertices \(x_0,x_1\), we obtain a natural map \(f_{x_0,x_1}:X(x_0,x_1)\to Y(f(x_0),f(x_1))\).  Indeed, it is for this reason that we call \(\psh{\Delta\wr\C}\) the category of \(\psh{\C}\)-enriched simplicial sets.
\subsection{Simplicial mapping objects and $A$-simplices}

\begin{defn} If \(S\) is a simplicial set equipped with a pair of vertices \((s_0,s_1):\Delta_0\coprod \Delta_0 \to X,\) we define \(S(s_0,s_1)\) to be the pullback of the diagram \[\Delta_0 \overset{(s_0,s_1)}{\to} S^{\partial \Delta_1} \leftarrow S^{\Delta_1},\] and we call it the \dfn{simplicial set of edges from \(s_0\) to \(s_1\)}.  This association is functorial in the category of bipointed simplicial sets and admits a right adjoint \(\Sigma\), the unreduced suspension functor, \[K\mapsto \left(\partial \Delta_1\to  K\times \Delta_1 \coprod_{K\times \partial \Delta_1} \Delta_0 \times \partial \Delta_1\right)\]
\end{defn}

\begin{lemma}For any \(\psh{\C}\)-enriched simplicial set \(X\) equipped with two vertices \((x_0,x_1)\), we may construct a simplicial presheaf \(\Map_X(x_0,x_1)\) on \(\C\), functorial in bipointed objects of \(\psh{\Delta\wr\C}\), such that \[X(x_0,x_1)=\Map_X(x_0,x_1)_0\] and \[\Hom(A,\Map_X(x_0,x_1))=\mathfrak{M}(A,X)(x_0,x_1).\]  Moreover, this functor arises from a cosimplicial enlargement of the functor \(A\mapsto \Delta_1[A]\).
\end{lemma}
\begin{proof}
We will show that \(\Sigma K [-]:\psh{\C}\to \psh{\Delta\wr \C}\) is a parametric left adjoint for any simplicial set \(K\).  It suffices to prove this when \(K\) is a simplex or empty, since \(\Sigma\) is well-known to be a parametric left adjoint.  The case when \(K\) is empty is clear, since \(\Sigma\emptyset=\Delta_0 \coprod \Delta_0\), and \((\Delta_0\coprod \Delta_0)[A]\) is just a coproduct of two vertices for every presheaf \(A\) on \(\C\).  

The case when \(K=\Delta_0\) is simply the functor \(\Delta_1[-]\), which is a parametric left adjoint by \eqref{leftadjointness}.  For \(K=\Delta_n\), we can decompose \(\Sigma K\) using the prism decomposition for the product \(\Delta_n\times \Delta_1\).  The prism decomposition presents \(\Delta_n\times \Delta_1\) as the colimit 
\[\varinjlim\left (\Delta_{n+1} \overset{\delta_n}{\leftarrow} \Delta_n \overset{\delta_n}{\rightarrow} \Delta_{n+1} \overset{\delta_{n-1}}{\leftarrow} \dots \overset{\delta_1}{\rightarrow} \Delta_{n+1}\right ). \]

When we take the pushout of the diagram \[\Delta_n\times \Delta_1 \leftarrow \Delta_n \times \partial \Delta_1\to \Delta_0 \times \partial \Delta_1,\] together with the prism decomposition, we find that \(\Sigma(\Delta_n)\) can be identified with the colimit of the diagram 
\[E_{n+1}^{n} \leftarrow E_{n}^{n-1} \rightarrow E_{n+1}^{n-1} \leftarrow E_n^{n-2} \to \dots \leftarrow E_n^1 \to E_{n+1}^1,\]
where \(E_n^i\) is the colimit of the diagram \[\Delta_0 \coprod \Delta_0 \leftarrow \Delta_{i-1} \coprod \Delta_{n-i} \hookrightarrow \Delta_n, \] where the map \(\Delta_{i-1}\hookrightarrow \Delta_n\) is the face spanned by the vertices \([0,\dots, i-1]\), and \(\Delta_{n-i}\hookrightarrow \Delta_n\) is the face spanned by the vertices \( [i,\dots,n]\).  Since \((-)[A]\) preserves colimits, it will suffice to show that \(E_n^i[-]\) is a parametric left adjoint for any pair \((n,i)\) such that \(1\leq i \leq n-1\).  

By the cocontinuity of \((-)[A]\), we may decompose \(E_n^i[-]\) as the pushout of the diagram \[\Delta_0\coprod \Delta_0 \leftarrow \Delta_{i-1} \coprod \Delta_{n-i}[-] \rightarrow \Delta_n[-].\]  However, it is clear from this construction that we may replace \(\Delta_{i-1} \coprod \Delta_{n-i}[-]\) by \[V[(i-1)+1+(n-i)](\ast,\dots,\ast,\emptyset,\ast,\dots,\ast)\] and \(\Delta_n[-]\) by \[V[(i-1)+1+(n-i)](\ast,\dots,\ast,-,\ast,\dots,\ast),\] since these are the parts of the functor that are killed in the pushout.  However, by \eqref{leftadjointness}, these functors are parametric left adjoints whose values on \(\emptyset\) are all exactly \(\partial \Delta_1=\Delta_1[\emptyset]\).  

It follows from this that we may define the aforementioned functor \[\Map:\overcat{\Delta_1[\emptyset]}{\psh{\Delta\wr \C}}\to \psh{\Delta\times \C}\] by the formula \[\Map_X(x_0,x_1)_n(c)=\Hom_{\overcat{\partial \Delta_1}{\psh{\Delta\wr \C}}}\left(\Sigma(\Delta_n)[c]_0^1, X_{x_0}^{x_1}\right),\] which is well-behaved since the functor \(\Sigma(K)[-]\) is a parametric left adjoint for every simplicial set \(K\).  We can see that \(\Hom(A,\Map_X(x_0,x_1))=\mathfrak{M}(A,X)(x_0,x_1)\) by unraveling the definitions.
\end{proof}

\subsection{The $\Delta\wr \C$-localizer $\mathsf{W}_{\Sc}$}
We begin with a small warning regarding notation:

\begin{note} Given a family of objects \(c=(c_1,\dots,c_n)\) of \(\C\), we will denote the presheaf \(Y_{\Delta\wr \C}([n](c_1,\dots, c_n))\) by \(\Delta_n[c]\).  Similarly, for a family of presheaves \((A_1,\dots,A_n)\) on \(\C\), we will denote \(V[n](A_1,dots,A_n)\) simply by \(\Delta_n[A]\).  We warn the reader that when \(A\) is simply a presheaf on \(\C\), this notation is used to mean \(V[n](A,\dots,A),\) but we are quite confident that the reader will be able to sort out which means which from context.  We just thought we'd let the reader know as a matter of courtesy.
\end{note}

\begin{defn}
Given a family \(c=(c_1,\dots,c_n)\) of objects of \(\C\), we define \dfn{Segal core} of the \(c\)-simplex \(\Delta_n[c]\) to be the \(\psh{\C}\)-enriched simplicial set 
\[\Sc_n[c]=\varinjlim\left(\Delta_1[c_1] \overset{\delta_0}{\leftarrow} \Delta_0 \overset{\delta_1}{\to} \dots \overset{\delta_0}{\leftarrow} \Delta_0 \overset{\delta_1}{\to} \Delta_1[c_n]\right ).\]
\end{defn}
\begin{defn} We define \(\mathsf{W}_{\Sc}\) to be the regular completion of the \(\Delta\wr\C\)-localizer generated by the class comprising the Segal core inclusions \(\Sc_n[c]\hookrightarrow \Delta_n[c]\) for any family of objects \(c=(c_1,\dots,c_n)\) of \(\C\). 
\end{defn}
\begin{lemma} The \(\psh{\C}\)-enriched simplicial set \(J=J[e]\), where \(J\) is the simplicial set obtained by taking the nerve of the strictly contractible groupoid \(G_2\) with two objects, is an injective object in \(\psh{\Delta\wr\C}\). 
\end{lemma}
\begin{proof}  The functor \(p:\Delta\wr \C \to \Delta=\Delta\wr \ast,\) induced by the terminal functor \(\C\to \ast,\) gives rise to an adjunction \[p_!:\psh{\Delta\wr\C}\leftrightarrows\psh{\Delta}:p^*.\]  We can see easily that \(p^*(X)=X[e]\) for any simplicial set \(X,\) because the functor \(p^*\) itself admits a right adjoint, which is \(\mathfrak{M}(e,-)\).  Then \(J[e]=p^*J=p^*\mathfrak{N}_\Delta(G_2),\) so it will suffice to show that \(p^*\mathfrak{N}_\Delta\) sends trivial fibrations in the natural model structure on \(\cat{Cat}\) to trivial fibrations of \(\psh{\C}\)-enriched simplicial sets.

However, this is equivalent to asking that the left adjoint of this functor sends monomorphisms of \(\psh{\C}\)-enriched simplicial sets to cofibrations between categories.  However, cofibrations in \(\cat{Cat}\) are just functors that induce injections on sets of objects.  We leave the easy proof of this fact to the reader.  
\end{proof}
\begin{cor}For every \(\psh{\C}\)-enriched simplicial set \(X\), the canonical map \(J\times X\to X\) is a trivial fibration, and in particular, belongs to \(\mathsf{W}_{\Sc}\).
\end{cor}
\begin{proof} Since \(J\) is an injective object, the map \(J\to e\) is a trivial fibration, which means that the map \(X\times J\to X\) is a trivial fibration as well, and therefore, it follows that \(X\times J\to X\) belongs to \(\mathsf{W}_\wr,\) since localizers contain all trivial fibrations.  
\end{proof}
\begin{cor}\label{ssetadjunction}The functor \((-)[e]:\psh{\Delta}\to \psh{\Delta\wr \C}\) is a left Quillen functor when \(\psh{\Delta}\) is equipped with the Joyal model structure and when \(\psh{\Delta\wr\C}\) is equipped with the Cisinski model structure generated by \(\mathsf{W}_{\Sc}\).  
\end{cor}
\begin{proof} Since the functor \((-)[e]=p^*\) admits an exceptional left adjoint, it necessarily preserves monomorphisms.  For this functor to preserve weak equivalences, we may equivalently show that the preimage of \(\mathsf{W}_\wr\) contains the Joyal weak equivalences. We can show that this is the case, then, by showing that the preimage is itself a \(\Delta\)-localizer containing the spine inclusions, which are known to generate the Joyal weak equivalences.  However, by \cite{cisinski-book}*{Proposition 1.4.20}, the preimage forms a \(\Delta\)-localizer provided that there exists some functorial cylinder \(\mathfrak{I}=(I,\partial^0,\partial^1,\sigma)\) of \(\Delta\) such that \(\sigma_X[e]:(I\otimes X)[e]\to X[e]\) is belongs to \(\mathsf{W}_\wr\) for every simplicial set \(X\).  

Since the functor \((-)[e]\) preserves products, again, since it admits a left adjoint, we see that the projection \((X\times J)[e]\to X[e]\) is exactly \(X[e]\times J[e]\to X[e]\), which belongs to \(W_\wr\) by the previous corollary.  This implies that the preimage of \(W_\wr\) indeed forms a \(\Delta\)-localizer, and this localizer clearly contains the spine inclusions, since these are mapped to Segal cores.
\end{proof}
We need to make use of a technical but straightforward lemma in order to obtain the upshot:
\begin{lemma}\label{scompcart} If \(\W\) is an accessible regular \(\C\)-localizer for a small category \(\C\), then \(\W\) is cartesian if and only if its simplicial completion\(\W_\Delta\) is cartesian.  
\end{lemma}
\begin{proof} Suppose \(\W\) is cartesian.  Then since the localizer is regular, its simplicial completion \(\W_\Delta\) is the class of weak equivalences obtained as the weak equivalences of the left Bousfield localization of the injective model structure on simplicial presheaves at the set \(S\times \Delta_0= \{s\times \Delta_0: s\in S\}\) for some set of maps \(S\) generating \(\W\).  Since the class of weak equivalences of the injective model structure is cartesian, it suffices to show that for any simplicial presheaf \(T\) on \(\C\) and any map \(s:A\to B\) in \(S\), the map \(s\times T: A\times T \to B\times T\) belongs to \(\W_\Delta\).  By regularity, \(T\) is the homotopy colimit of its category of elements, so \(s\times T\) is the homotopy colimit of maps of the form \(s\times (c \times \Delta_i)\), where the \(c\times \Delta_i\to T\) is a section of \(T\) for some object \((c,i)\) of \(\C\times \Delta\).  Then we have that \[s\times (c \times \Delta_i)
= (s\times c)\times \Delta_i,\] and since \(s\times c\) belongs to \(\W\) for every \(c\) in \(\C\), all of these maps are objectwise \(\W\)-equivalences and therefore belong to \(\W_\Delta\).  Since \(s\times T\) is the homotopy colimit of a diagram of weak equivalences, it is itself a weak equivalence and therefore, \(\W_\Delta\) is cartesian.  

The converse follows immediately from \cite{cisinski-book}*{Proposition 2.3.37}.  \end{proof}
\begin{lemma}\label{weakcats} The regular \(\Delta\wr \C\)-localizer \(\mathsf{W}_{\Sc}\) generated by the Segal cores is accessible and cartesian.  Moreover, it is strongly generated by the set of maps comprising the Segal core inclusions and the map \(j:J\to e\).  
\end{lemma}
\begin{proof} 
The first assertion is proven in two separate parts, since by \eqref{boostgenerators}, we note that \(\W_{\Sc}\) is strongly generated by the class \(\Sc \cup \operatorname{cart}(\{j\})\), where \(\Sc\) denotes the set of Segal cores.  Then we first show that \(\Sc\times \Delta_0\) generates a cartesian \(\Delta \wr \C \times \Delta\)-localizer.     
 
This is exactly the content of \cite{rezk-theta-n-spaces}*{Theorem 6.6} because we are looking at the regular completion, which means that we are Bousfield localizing the class of discrete Segal cores over the injective model structure. 

However, the reader should beware that the proof depends on \cite{rezk-theta-n-spaces}*{Proposition 6.4}, which was left uncorrected in the most recent revision of the paper.  The proof stated there is based on an incorrect statement from the published revision, and the author had forgotten to update it in the correction.  However, the proof of \eqref{coversweak} later in this paper can easily be modified to give a correct proof of that assertion. 

That the cartesian property holds for the whole simplicial completion is a corollary of \eqref{scompcart}, since this implies that \(\operatorname{cart}(\{j\})\times \Delta_0\) generates a cartesian \(\Delta \wr \C \times \Delta\)-localizer, and by \cite{cisinski-book}*{Corollary 1.4.19b}, \[\W(\Sc\times \Delta_0 \cup \mathsf{R}(\W_\infty))\cup \W(\operatorname{cart}(\{j\})\times \Delta_0 \cup \mathsf{R}(\W_\infty))\] generates a cartesian \(\Delta \wr \C \times \Delta\)-localizer, since each of the two parts generate cartesian \(\Delta \wr \C \times \Delta\)-localizers.  This implies by \eqref{scompcart} that \(\W_{\Sc}\) is indeed cartesian.

The second claim the real content of \cite{rezk-theta-n-spaces}*{Proposition 7.21}, and we refer the reader to the proof given there.
\end{proof}

\begin{lemma}The functor \(\Sigma(-)[A]:\Delta \to \psh{\Delta\wr\C}\) is a functorial cosimplicial resolution for the functor \(\Delta_1[A]\) \(\C\) associated with the localizer \(\mathsf{W}_{\Sc}\).  Moreover, this same cosimplicial resolution is also a cosimplicial resolution for \((\Delta_1)[A]\) viewed as an object in the coslice category under \(\Delta_1[\emptyset]\).
\end{lemma}
\begin{proof} First, we can see that the map \(\Delta[1]\to \Sigma\Delta[n]\) is inner anodyne as follows: First, we may form the pushout product of the spine inclusion \(\iota_n:\Sp[n]\hookrightarrow \Delta[n]\) with the monomorphism \(b:\partial\Delta[1]\hookrightarrow \Delta[1]\).  Since inner anodyne maps are closed under pushout-products, this gives us an inner-anodyne map \[\iota_n \wedge b: \Delta[1]\times \Sp[n] \cup \partial\Delta[1]\times \Delta[n]\hookrightarrow \Delta[n]\times \Delta[1].\]  

However, the source of this map admits another canonical map induced by the commutativity of the square under projection to the suspension of the spine, which is bipointed.  That is to say, we are looking at the canonical map \[m:\Delta[1]\times \Sp[n] \cup \partial\Delta[1]\times \Delta[n] \to \Sigma \Sp[n].\]  Pushing out the pushout-product map along this map \(m\), we obtain an inner anodyne map from the suspension of the spine to the suspension of the \(n\)-simplex.  To see that the inclusion of any nontrivial edge into the suspension of the spine is inner anodyne, suppose we're given a lifting diagram with the inclusion of a nontrivial edge into the suspension of the spine on the left together with an inner fibration against which we must find a lift. However, since the suspension of the spine is just a finite-length family of 2-disks glued together along their opposite edges, we may continually extend the original edge along degenerate edges in front or behind and thereby find a lift of each disk by induction.

Then in particular, the maps \(\Sigma\Delta[n]\to \Delta[1]\) are retracts of the inner anodyne inclusions \[\Delta[1]\hookrightarrow \Sigma\Sp[n]\hookrightarrow \Sigma \Delta[n]\],  and therefore the functors \((-)[A]\) must send them to weak equivalences, since those functors necessarily preserve inner anodyne maps, as they send spine inclusions to Segal core inclusions.   

It suffices then to show that for any presheaf \(A\) on \(\C\), the cosimplicial object defined by the functor \(\Sigma(-)[A]\) is Reedy cofibrant both as an ordinary resolution and as a cosimplicial resolution in the coslice under \(\Delta_1[\emptyset]\).  The second case is immediate, since the coslice version of the cosimplicial resolution preserves monomorphisms and colimits.  Then we consider the other case.

However, this case is similarly trivial because the functor preserves connected colimits of simplicial sets and monomorphisms, so, in particular, the image of the boundary of the \(n\)-simplex injects into the \(n\)th component, which gives that the latching map is a monomorphism, and therefore that the object is a Reedy cofibrant cosimplicial object.  
\end{proof}
\begin{cor}For any \(\mathsf{W}_{\Sc}\)-fibrant \(\psh{\C}\)-enriched simplicial set \(X\) equipped with two vertices \(x_0,x_1\), the simplicial set \(\Map_X(x_0,x_1)(A)\) models the homotopy function complex whose set of connected components is \(\psh{\C}[\mathsf{W}_{\Sc}^{-1}]((\Delta_1[A],0,1),(X,x_0,x_1))\). Similarly, when \(X\) is \(\mathsf{W}_{\Sc}\)-fibrant, the simplicial set \(E(A,X)_n=\Hom(\Sigma(\Delta_n)[A],X)\) gives a model for the homotopy function complex \(H_{\mathsf{W}_{\Sc}}(\Delta_1[A],X)\).  
\end{cor}
\begin{proof}This follows immediately from the preceding lemma, since we have merely constructed homotopy function complexes from the given resolutions.  
\end{proof}

\subsection{The $\Delta\wr \C$-localizer $\mathsf{W}_{\wr}$}\label{weakenrichment}
We fix an accessible cartesian regular \(\C\)-localizer \(\W\).  We will give the definition of the \(\Delta\wr C\)-localizer \(\W_\wr\), and using \cite{rezk-theta-n-spaces}*{Theorem 8.1}, we will show that it is cartesian.  Moreover, we will show that if \(S\) is a class of strong generators for \(\W\), then \(\W_{\wr}\) is strongly generated (over \(\W_\Sc\)) by \(\Delta_1[S]\), the class comprising those maps \(\Delta_1[f]:\Delta_1[A]\to \Delta_1[B]\) such that \(f:A\to B\) belongs to \(S\).    
\begin{defn}We define a \dfn{suspended \(\W\)-equivalence} to be a map of the form \(\Delta_1[f]:\Delta_1[A]\to \Delta_1[B]\), where \(f:A\to B\) belongs to \(\W\).  
\end{defn}
\begin{defn} We define \(\mathsf{W}_\wr\) to be the regular completion of the \(\Delta\wr\C\)-localizer generated by the class comprising:
\begin{enumerate}
\item [(i)] The suspended \(\mathsf{W}\)-equivalences.
\item [(ii)] The Segal core inclusions \(\Sc_n[c]\hookrightarrow \Delta_n[c]\) for any family of objects \(c=(c_1,\dots,c_n)\) of \(\C\).
\end{enumerate}     
\end{defn}
Before we begin, we first fix some notation, to avoid confusion.  We will denote the \emph{regular} \(\Delta\wr \C\times \C\)-localizer generated by a class of maps \(S\times \Delta_0\) with \(S\) a class of maps in \(\psh{\Delta\wr \C\times \Delta}\) by \(\sW(S\cup \W_\infty)\), which coincides with the localizer \(\W(S\cup \mathsf{R}(\W_\infty))\).  If \(S\) is a class of maps in \(\psh{\Delta\wr \C}\), we will, by abuse of notation, let \(\sW(S)=\sW(S\times \Delta_0)\).  
\begin{thm}[\cite{rezk-theta-n-spaces}*{Proposition 8.5}] The regular \(\Delta\wr \C\times \Delta\)-localizer \(\sW(\Sc\cup \Delta_1[\W])\) is cartesian.  Moreover, if \(S\) is a class of strong regular generators for \(\W\), then \(\sW(\Sc\cup \Delta_1[\W])=\sW(\Sc\cup \Delta_1[S])\)
\end{thm}
\begin{proof}See \cite{rezk-theta-n-spaces}*{Proposition 8.2-8.5}.
\end{proof}
This immediately gives us the corollary
\begin{cor}The localizer \(\W_\wr\) is cartesian, and if \(S\) is a class of strong regular generators of \(\W\), then \(\Delta_1[S]\cup \Sc \cup \{j\}\) is a class of strong regular generators for \(\W_\wr\).   
\end{cor}

\section{The theory of $\G$-extensions and strict $\omega$-categories}

This section is mainly meant to be a quick review of the main results in the second and third chapters of \cite{ara-thesis}, and the author makes no claims of originality in this section. 

\subsection{The globe category $\G$}

A good deal of this section is taken straight from the first chapter of Dimitri Ara's thesis, \cite{ara-thesis}.

Let \(\mathbb{G}_n\) denote the category presented as the free category on 
\begin{equation*}
\begin{tikzpicture}
\matrix (m) [matrix of math nodes, row sep=3em,
column sep=3em, text height=1.5ex, text depth=0.25ex]
{ D_0 & D_1 & \dots &D_{n-1}& D_n \\};
\path[->, font=\scriptsize]
(m-1-1) edge[transform canvas={yshift=0.2em}]   node[auto]{\(\scriptstyle \sigma_1\)} (m-1-2)
(m-1-1) edge[transform canvas={yshift=-0.2em}]   node[auto, swap]{\(\scriptstyle \tau_1\)} (m-1-2)
(m-1-2) edge[transform canvas={yshift=0.2em}]  node[auto]{\(\scriptstyle \sigma_2\)}(m-1-3)
(m-1-2)	edge[transform canvas={yshift=-0.2em}]  node[auto, swap]{\(\scriptstyle \tau_2\)} (m-1-3)
(m-1-3) edge[transform canvas={yshift=0.2em}]  node[auto]{\(\scriptstyle \sigma_{n-1}\)}(m-1-4)
(m-1-3)	edge[transform canvas={yshift=-0.2em}] node[auto, swap]{\(\scriptstyle \tau_{n-1}\)} (m-1-4)
(m-1-4) edge[transform canvas={yshift=0.2em}]  node[auto]{\(\scriptstyle \sigma_{n}\)}(m-1-5)
(m-1-4)	edge[transform canvas={yshift=-0.2em}] node[auto, swap]{\(\scriptstyle \tau_{n}\)} (m-1-5);
\end{tikzpicture}
\end{equation*}
modulo the the coglobular relations, \[\sigma_{i+1}\sigma_i=\tau_{i+1}\sigma_i \qquad \text{and} \qquad \tau_{i+1}\tau_i=\sigma_{i+1}\tau_i\] for \(1\leq i\leq n\).  There is an obvious inclusion map \(\mathbb{G}_n\hookrightarrow \mathbb{G}_{n+1}\) for each \(n \in \mathbf{N}\).  This defines a directed system, and we denote its colimit in \(\cat{Cat}\) by \(\mathbb{G}\).  

For integers \(0\leq i\leq j\), we define maps \(D_i\to D_j\) in \(\mathbb{G}\): 

\[\sigma^j_i=\sigma_j\hdots\sigma_{i+1}\qquad \text{and}\qquad\tau^j_i=\tau_j\hdots\tau_{i+1}\]

It follows by induction and the coglobular relations that given \(D_n, D_m\in \mathbb{G}\), we have that 
\[\Hom_\mathbb{G}(D_n,D_m)=
\begin{cases}\{\sigma_n^m,\tau_n^m\} & \text{if \(n<m\)}\\
\{\id_{D_n}\} & \text{if \(n=m\)}\\
\emptyset & \text{otherwise}\end{cases} \]

For any presheaf in \(X\in \operatorname{Ob}\hat{\G}=\cat{Cat}(\G^\op,\cat{Set})\), by abuse of notation, we let \(s_n=X_{\sigma_n}\) and \(t_n=X_{\tau_n}\).

\subsection{Globular patterns and $\G$-extensions}

For \(k\geq 2\), we define the category \(I_k\) to be the category associated with the ordered set \(\{(i,j): 0\leq i \leq 1 \wedge 0\leq j\leq k \wedge (i,j)\neq (0,k)\}\) ordered by the relation that \((i,j)\leq (i',j')\) if and only if \(i'-i=1\) and \(0\leq j'-j\leq 1\).  When \(k=1\), we let \(I_k=\ast\). 

\begin{defn}
A functor \(\eta:I_k\to \mathbb{G}\) for \(k\geq 1\) is called a \dfn{globular pattern} when the following conditions are satisfied:

\begin{enumerate}
\item [(i)] Every morphism of the form \(\alpha:(0,j)\to (1,j)\) in \(I_k\), the map \(\eta(\alpha)=\sigma_n^m\) for some \(m>n\geq 0\)
\item [(ii)] Every morphism of the form \(\beta:(0,j)\to (1,j+1)\) in \(I_k\), the map \(\eta(\beta)=\tau_n^m\) for some \(m>n\geq 0\)
\end{enumerate}

If \((C,F:\G\to C)\) is a category under \(\G\), a functor \(\eta:I_k\to C\) for some \(k\geq 1\) such that \(\eta\) factors as \(\eta_0^\ast F=F\circ \eta_0\) for some globular pattern \(\eta_0:I_k\to \G\) is called a \dfn{globular pattern in \((C,F)\)}.

If \((C,F:\G^\op\to C)\) is a category under \(\G^\op\), we define a \dfn{coglobular pattern in \((C,F)\)} to be a functor \(\eta:I_k^\op \to C\) for some \(k\geq 1\) such that the corresponding functor \(\eta^{op}:I_k\to C^\op\) is a globular pattern in \((C^\op,F^\op:\G\to C^\op)\).

We define globular sums (resp. globular products) in a category \((C,F:\G\to C)\) under \(\G\) (resp. in a category \((C,F:\G^\op\to C)\) under \(\G^\op\)), to be colimits (resp. limits) of globular patterns (resp. coglobular pattern) \(\eta\) in \((C,F)\).  
\end{defn} 

\begin{defn}
We say that a category \((C,F:\G\to C)\) under \(\G\) (resp. in a category \((C,F:\G^\op\to C)\) under \(\G^\op\)) to be a \dfn{globular \(\G\)-extension} (resp. \dfn{ globular \(\G\)-coextension}) if it contains all globular sums (resp. globular products).  A morphism of \(\G\)-extensions is a functor under \(\G\) that preserves all  globular sums.  Unless otherwise noted, we will simply refer to these as \(\G\)-extensions and \(\G\)-coextensions respectively.
\end{defn}
 
\begin{defn}
Given a \(\G\)-extension \((C,F)\) and a category \(D\), we define a \dfn{\(D\)-valued \(C\)-model} to be a functor \(G:C^\op\to D\) such that \((D,G\circ F^\op:\G^\op\to D)\) is a \(\G\)-coextension and such that the functor \(G^{op}:C\to D^\op\) is a morphism of \(\G\)-extensions.  We define the category \(\Mod(C,D)\) to be the full subcategory of \(D^{C^\op}\) spanned by the \(D\)-valued \(C\)-models.  

By abuse of notation, for any category \(D\) we denote the full subcategory of \(D^{\G^\op}\) spanned by the \(\G\)-coextensions by \(\Mod(\G,D)\).  We call the objects of this category \dfn{globular sets taking values in \(D\)}.
\end{defn}

\subsection{The initial $\G$-extension $\Theta_0$}

\begin{prop}
There exists a unique \(\G\)-extension \((\Theta_0,\iota_0\:\G\to \Theta_0)\) such that the induced transformation 
\[\iota_0^*:\Mod(\Theta_0,\cdot)\to \Mod(\G,\cdot)\]
is an equivalence of 2-functors.  Moreover, for any \(\G\)-extension \(D,F\), there exists a unique (up to isomorphism) factorization of the structure map \(F\) as the composite of the map \(\iota_0:\G\to \Theta_0\) with some morphism of \(\G\)-extensions \(F_0:\Theta_0\to D\) .  
\end{prop}
\begin{proof}
We take \(\Theta_0\) to be the full subcategory of \(\widehat{\G}\) spanned by the globular sums.  This gives a \(\G\)-extension because \(\Theta_0\subseteq \widehat{\G}\) contains the image of the Yoneda embedding, which gives us a factorization \(\G\hookrightarrow \widehat{\G}=\G\overset{\iota_0}{\hookrightarrow} \Theta_0\overset{\gamma_0}{\hookrightarrow}\widehat{\G}\).   

Then we would like to construct an inverse for the transformation \(\iota_0^\ast\). Let \(\cat{Comp}\) be the 2-subcategory of \(\cat{Cat}\) spanned by the complete categories with limit-preserving functors between them.  By the universal property of the co-Yoneda embedding, we have that for any complete category \(B\), \(\cat{Cat}(\G^\op,B)\simeq \cat{Comp}(\widehat{\G}^{op},B)\), naturally in \(B\). Also, since every complete category \(B\) necessarily contains all globular products for any functor \(\G^\op\to B\), and since every continuous functor \(X:\widehat{G}^{op}\to B\) necessarily preserves all globular products \emph{a fortiori}, we have an embedding \(\cat{Comp}(\widehat{\G}^{op},B)\hookrightarrow \Mod(\widehat{\G},B)\)  

Then we see that we have a chain of transformations natural in \(C\), the composite of which we will call \(F_C\),
\[\Mod(\G,C)\overset{\iota}{\hookrightarrow} \cat{Cat}(\G^\op,\widehat{C}) \simeq \cat{Comp}(\widehat{G}^\op,\widehat{C})\hookrightarrow\Mod(\widehat{\G},\widehat{C})\overset{\gamma_0^\ast}{\to} \Mod(\Theta_0,\widehat{C}),\]
but for each globular set \(X\) taking values in \(C\), the object \(FX:\Theta_0\to \widehat{C}\) necessarily factors uniquely (up to specified isomorphism) through the Yoneda embedding, since the category \(C,i_0^\ast FX:\G\to C\) under \(\G\) contains all globular products, and the Yoneda embedding \(C\hookrightarrow \widehat{C}\) necessarily preserves them.  

This means, in particular, that \(F_C\) factors through the inclusion \[\Mod(\Theta_0,C)\hookrightarrow \Mod(\Theta_0,\widehat{C}).\]  We let \(H_C\) denote the factor of \(F_C\) going from \(\Mod(\G,C)\) to \(\Mod(\Theta_0,C)\) (naturally in \(C\).  Suppressing the \(C\), we see that \(H\) clearly inverse to \(\iota_0^\ast\) by way of the factorization.   The fact that \(\Theta_0\) is initial in the category of \(\G\)-extensions follows from the earlier claim by letting \(C=D^\op\).
\end{proof}

We recall a proposition of Ross Street in \cite{street}:

\begin{prop} There exists an order structure \(\bltri_A\) on \(\overcat{\G}{A}\) for any presheaf \(A\) on \(\G\) such that maps \(A\to B\) of presheaves induce order preserving maps \(\overcat{\G}{A}\to \overcat{\G}{B}\) and such that the associated ordered sets \((\overcat{\G}{X}, \bltri_X)\) of the objects \(X\) of \(\Theta_0\subset \widehat{\G}\) are finite and linearly ordered.
\end{prop}
\begin{proof}
We construct a functorial order structure on \(\El(A)=\operatorname{Ob}(\overcat{\G}{A})\) for a presheaf \(A\) on \(\G\) following Street in \cite{street} by taking the reflexive transitive closure \(\blacktriangleleft_X\) of the relation\(\prec_X\) defined such that given a map \(\alpha:D_n\to X\), 
\begin{equation}\label{streetord} \alpha \prec_X \beta \qquad \text{if and only if \(\alpha=s_{n+1}(\beta)\) or \(t_n(\alpha)=\beta\)}.\end{equation}\
To show that this order structure is functorial, it suffices to show that the function \(\El(f):\El(A)\to \El(B))\) preserves the order structure generated by the relation above.  However, this follows from the definition \eqref{streetord} and the commutativity of the induced functor with the source and target maps. 

We will show that given any globular pattern \(H:I_k \to \G\) with colimit \(X\) in \(\widehat{\G}\), the ordered set \((\El(X),\bltri_X)\) is linearly ordered.  We proceed by induction as follows: Assume that for every globular pattern \(I_j\to \G\) with \(j<k\), the claim holds for \(colim I_j\).  Then we note that we may decompose \(X\cong Y\coprod_{D_{i'_{k-1}}} D_{i_k} \), where \(Y\) is the colimit of the globular pattern \(\iota_0^*H:I_{k-1}\to \G\), where \(\iota_0:I_{k-1}\to I_k\) is the obvious inclusion on the first \(k-1\) components.  Note that this gives a canonical factorization \(D_{i'_{k-1}}\to Y\) as \(D_{i'_{k-1}} \to D_{i_{k-1}} \to Y\), where the first map is \(\sigma^{i_{k-1}}_{i'_{k-1}}\).  

Then if \(\gamma:D_{i_\gamma}\to X, \lambda: D_{i_\lambda}\to X\) are both maps factoring through either \(\alpha_k:D_{i_k}\to X\) or \(\alpha_Y:Y\to X\), then \(\gamma \bltri_X \lambda\) reduces to \(\gamma \bltri_{D_{i_{k}}}\lambda\) or \(\gamma\bltri_Y \lambda\).  So without loss of generality, since every element of \(\El(X)\) factors through at least one globular summand, we may assume that \(\gamma\) belongs to \(\El(D_{i_k}) - \im(\El(\tau^{i_k}_{i'_{k-1}}))\) and that \(\lambda\) belongs to \(\El(Y)-\im(\El(\sigma^{i_{k-1}}_{i'_{k-1}}))\).  

First, notice that for \(j\leq i'_{k-1}\), we have that \(s^{i'_{k-1}}_j(\alpha'_{k-1})=s^{i_{k-1}}_j(\alpha_{k-1})\) and that \(t^{i'_{k-1}}_j(\alpha'_1=t^{i_k}(\alpha_Y)\).  

Since \(\gamma\) belongs to \(\El(D_{i_k}) - \im(\El(\tau^{i_k}_{i'_{k-1}}))\), we have that \(\gamma\) lives in the strict \(\bltri_{D_{i_k}}\)-interval \((s^{i_k}_{i'_{k-1}}(\alpha_k),t^{i_k}_{i'_{k-1}}(\alpha_k))_{\bltri_{D_{i_1}}})\).  Similarly, since \(\lambda\) belongs to \(\El(Y)-im(\El(\sigma^{i_{k-1}}_{i'_{k-1}}))\), this in particular implies that \(s^{i_{k-1}}_{i'_{k-1}}(\alpha_{k-1}) \bltri_Y \lambda\).  

Since \[\gamma \bltri_{D_{i_k}} t^{i_k}_{i'_{k-1}}(\alpha_k),\] \[t^{i_k}_{i'_{k-1}}(\alpha_k)=\alpha'_{k-1}= s^{i_{k-1}}_{i'_{k-1}}(\alpha_{k-1}),\] and also \[s^{i_{k-1}}_{i'_{k-1}}(\alpha_{k-1})\bltri_Y \lambda,\] we have that \(\gamma \bltri \lambda\).  It is clear that antisymmetry holds, since in a case such as the one above, there is exactly one traversable path between components, and when they are in the same component, antisymmetry is inherited from the lower-order cases by induction. 
\end{proof} 

\begin{prop}
The category \(\Theta_0\) is the full subcategory of \(\psh{\G}\) spanned by the objects \(X\) such that the ordered set \[(\overcat{\G}{X}, \bltri_X)\] is finite and linearly ordered.
\end{prop}
\begin{proof}
Since \(\overcat{\G}{X}\) is finite, let \(i_1,\dots i_n\) be the heights corresponding to the elements of maximal height and ordered as a subset of the total order under \(\bltri_X\).  (Finish proof later)
\end{proof}

\begin{prop}
Every morphism \(f:A\to B\)  in \(\Theta_0\) is monic.
\end{prop}
\begin{proof} It suffices to show that the natural transformation between presheaves on \(\G\) associated with the map \(f\) is objectwise injective.  To see this, we give \(A_n\) the structure of a linear digraph, where \(\alpha \prec_n \beta\) if and only if \(\alpha \bltri 
\beta\) and \((\alpha,\beta)_\bltri\cap A_n=\emptyset\), where \((\alpha,\beta)\) denotes the strict open interval between \(\alpha\) and \(\beta\). We see that \(f_n\) preserves \(prec_n\), but \(prec_n\) is irreflexive, so in particular \(f_n\) is injective.    

\end{proof}

\subsection{The globular envelope of a category under $\Theta_0$}

Let \((C,D_C)\) be a \(\G\)-extension.  Then we say that a functor \(F:C\to E\) is a \dfn{globular \((C,D_C)\)-extension} if \((D,D_C^\ast(F))\) is a \(\G\)-extension and \(F\) is a morphism of \(\G\)-extensions.  We define the \dfn{category of globular \((C,D_C)\)-extensions}, denoted \((C,D_C)\cat{-Ext}\) to be the category whose objects are globular \((C,D_C)\)-extensions and whose arrows are morphisms of \(\G\)-extensions under \((C,D_C)\).  Unless otherwise noted, we will abuse notation and simply denote this category simply by \(C\cat{-Ext}\), with its objects similarly called \(C\)-extensions.  

\begin{prop} The category of \(\Theta_0\)-extensions is equivalent to the category of \(\G\)-extensions.  \end{prop}
\begin{proof} Immediate from the definitions.
\end{proof}

\begin{prop} If \((C,D_C)\) is a \(\G\)-extension, any functor under \(C\) between \(C\)-extensions is a morphism of \(C\)-extensions.
\end{prop}
\begin{proof} Let \[H:(X,F_X)\to(Y,F_Y)\] be a functor under \(C\) between \(C\)-extensions, and let 
\begin{align*}D_X&=D_C^\ast(F_X) \intertext{and} D_Y&=D_C^\ast(F_Y).\end{align*}  Then any globular sum in \((X,D_X)\) is the image under \(F_X\) of a globular sum in \((C,D_C)\). Then since \(HF_X=F_Y\), and \(F_Y\) sends globular sums in \((C,D_C)\) to globular sums in \((Y,D_Y)\), it follows that any globular sum in \((X,D_X)\) must map under \(H\) to a globular sum in \((Y,D_Y)\).  Therefore, \(H\) is a morphism of \(\G\)-extensions and belongs to \((C,D_C)\cat{-Ext}\).  
\end{proof}
This immediately yields the corollary:
\begin{cor} For any \(\G\) extension \((C,D_C)\), the category \(C\cat{-Ext}\) is a full subcategory of \(\overcat{C}{\cat{Cat}}\)\end{cor}.

Until the close of this subsection, we denote the forgetful functor \(\Theta_0\cat{-Ext} \to \overcat{\Theta_0}{\cat{Cat}}\) by the letter \(U\).

\begin{defn}Given a category \(C\) equipped with a functor \(\Theta_0\to C\), we say that a functor \(C\to U(C')\) under \(\Theta_0\) \dfn{exhibits \(C'\) as a globular envelope of \(C\)} if the following property holds: 

Given any solid arrow diagram
\begin{equation*}
\begin{tikzpicture}
\matrix (b) [matrix of math nodes, row sep=3em,
column sep=1.5em, text height=1.5ex, text depth=0.25ex]
{   & C & \\
  U(C') & & U(D) \\};
\path[->, font=\scriptsize]
(b-1-2) edge (b-2-1)
        edge (b-2-3)
(b-2-1) edge [dotted] node [auto,swap] {\(\scriptstyle U(f)\)} (b-2-3);
\end{tikzpicture},
\end{equation*}
there exists a unique arrow \(f:C'\to D\) in \(\Theta_0\cat{-Ext}\) such that \(U(f)\) gives the desired dotted arrow.  In such a situation, we will call \(C'\) a \dfn{globular envelope} for \(C\).  It is clear from the definition that any two globular envelopes for \(C\) are unique up to unique isomorphism.  
\end{defn}

\begin{thm}[\cite{ara-thesis}*{2.6}] Every small category under \(\Theta_0\) admits a globular envelope.
\end{thm}
\begin{proof} See \cite{ara-thesis}*{2.6}.   
\end{proof}

\subsection{Categorical $\G$-extensions}

A categorical \(\G\)-extension should, to a first approximation, be a \(\G\)-extension \(F:\G\to C\) together with the data of co-composition and co-degeneracy morphisms endowing the corresponding globular classified by the \(\G\)-coextension \[F^\op:\G^\op \to C^\op\] with the structure of a strict \(\omega\)-category internal to \(C^{op}\).  We write out what this means explicitly:

\begin{defn}
For \(i\geq j\geq 0\)
A \dfn{precategorical \(\G\)-extension} is specified by the following data:
\begin{enumerate}
\item [(i)] A functor \(D:\G \to C\) equipping \(C\) with the structure of a \(\G\)-extension.  We write \(D(D_n)\) simply as \(D_n\), and for \(f\in \G\), simply write \(D(f)\) as \(f\).  
\item [(ii)] For each \(i>j\geq 0\), a morphism \(\nabla^i_j:D_i \to D_i \coprod_{D_j} D_i\)
\item [(iii)] For each \(i\geq 0\), a morphism \(\kappa_i:D_{i+1}\to D_i\)
\end{enumerate}
satisfying the following axioms:
\begin{enumerate}
\item [(PC1)] For each \(i> 0\), we have that \[\kappa_i \sigma_{i+1}=\id_{D_i} \qquad \text{and} \qquad   \kappa_i \tau_{i+1}=\id_{D_i}\]
\item [(PC2)] For each \(i>j\geq 0\), let \(\varepsilon_1\) and \(\varepsilon_2\) denote the two canonical maps \(D_i\to D_i\coprod_{D_j} D_i\), we have that: 
\begin{align*}
\nabla^i_j \sigma_i = &\begin{cases}\varepsilon_2\sigma_i & \text{if \(j=i-1\)}\\ (\sigma_i\coprod_{D_j} \sigma_i)\nabla^{i-1}_j &\text{otherwise}\end{cases}\\
\intertext{and}\\
\nabla^i_j \tau_i = &\begin{cases}\varepsilon_1\tau_i & \text{if \(j=i-1\)}\\ (\tau_i\coprod_{D_j} \tau_i)\nabla^{i-1}_j & \text{otherwise}\end{cases}
\end{align*}
\end{enumerate}

\end{defn}

In keeping with Ara's treatment, we will fix the notations \[\kappa^j_i=\kappa_j\hdots\kappa_{i-2}\kappa_{i-1}\qquad\text{and}\qquad \nabla_k=\nabla^k_{k-1},\]
for \(i\geq j\geq 0\) and \(k>0\), respectively.  

\begin{defn} Using the same notation as above, we say that a precategorical \(\G\)-extension is a \dfn{categorical \(\G\)-extension} if it satisfies the following axioms:
\begin{enumerate}
\item[(CC1)] Associativity:\\
For \(i>j\geq 0\), the diagram
\begin{equation*}
\begin{tikzpicture}
\matrix (ass) [matrix of math nodes, row sep=4em,
column sep=4em, text height=1.5ex, text depth=0.25ex]
{ D_i & D_i\coprod_{D_j}D_i \\
   D_i\coprod_{D_j}D_i &  D_i\coprod_{D_j}D_i\coprod_{D_j}D_i \\};
\path[->, font=\scriptsize]
(ass-1-1) edge node[auto]{\(\scriptstyle \nabla^i_j\)} (ass-1-2)
          edge node[auto,swap]{\(\scriptstyle \nabla^i_j\)} (ass-2-1)
(ass-2-1) edge node[auto]{\(\scriptstyle \id_{D_i}\coprod_{D_j}\nabla^i_j\)} (ass-2-2)
(ass-1-2) edge node[auto,swap]{\(\scriptstyle \nabla^i_j\coprod_{D_j} \id_{D_i}\)} (ass-2-2);
\end{tikzpicture}
\end{equation*}
commutes.
\item[(CC2)]Strict interchange:\\
For \(i>j>k\geq 0\), the diagram
\begin{equation*}
\begin{tikzpicture}[bij/.style={above,sloped,inner sep=0.5pt}]
\matrix (ich) [matrix of math nodes, row sep=2.5em,
column sep=1.3em, text height=1.5ex, text depth=0.25ex]
{
                  & D_i & \\
D_i\coprod_{D_k}D_i & & D_i\coprod_{D_j}D_i \\
(D_i\coprod_{D_j}D_i)\coprod_{D_k} (D_i\coprod_{D_j}D_i) & & (D_i\coprod_{D_k}D_i)\coprod_{D_j\coprod_{D_k}D_j} (D_i\coprod_{D_k}D_i) \\};
\path[->, font=\scriptsize]
(ich-1-2.240) edge node[auto,swap]{\(\scriptstyle \nabla^i_k\)} (ich-2-1.60)
(ich-1-2.300) edge node[auto]{\(\scriptstyle \nabla^i_j\)} (ich-2-3.120)
(ich-2-1) edge node[auto]{\(\scriptstyle \nabla^i_j\coprod_{D_k}\nabla^i_j\)} (ich-3-1)
(ich-2-3) edge node[auto]{\(\scriptstyle \nabla^i_k\coprod_{\nabla^j_k}\nabla^i_k\)} (ich-3-3)
(ich-3-1) edge node[bij] {\(\sim\)}
							 node[below] {\(\scriptstyle \Phi\)} (ich-3-3);
\end{tikzpicture}
\end{equation*}
commutes, where \(\Phi\) is the unique isomorphism between the objects \[(D_i\coprod_{D_j}D_i)\coprod_{D_k}(D_i\coprod_{D_j}D_i)\] and \[(D_i\coprod_{D_k}D_i)\coprod_{D_j\coprod_{D_k}D_j} (D_i\coprod_{D_k}D_i)\] viewed as cones on the diagram
\begin{equation*}
\begin{tikzpicture}
\matrix (ijk) [matrix of math nodes, row sep=2em,
column sep=2em, text height=1.5ex, text depth=0.25ex]
{
D_i& D_k & D_i \\
D_j & D_k & D_j \\
D_i & D_k & D_j \\};
\path[->, font=\scriptsize]
(ijk-1-2) edge node[auto,swap]{\(\scriptstyle \sigma_k^i\)} (ijk-1-1)
					edge node[auto]{\(\scriptstyle \tau_k^i\)} (ijk-1-3)
(ijk-2-1) edge node[auto]{\(\scriptstyle \sigma_j^i\)} (ijk-1-1)
					edge node[auto,swap]{\(\scriptstyle \tau_j^i\)} (ijk-3-1)
(ijk-3-2) edge node[auto]{\(\scriptstyle \sigma_k^i\)} (ijk-3-1)
					edge node[auto,swap]{\(\scriptstyle \tau_k^i\)} (ijk-3-3)
(ijk-2-3) edge node[auto,swap]{\(\scriptstyle \sigma_j^i\)} (ijk-1-3)
					edge node[auto]{\(\scriptstyle \tau_j^i\)} (ijk-3-3)
(ijk-2-2) edge node[auto,swap]{\(\scriptstyle \sigma_k^j\)} (ijk-2-1)
					edge node[auto]{\(\scriptstyle \tau_k^j\)} (ijk-2-3)
(ijk-2-2) edge[-,double equal sign distance] node[auto,swap]{\(\scriptstyle \id_{D_k}\)} (ijk-1-2)
					edge[-,double equal sign distance] node[auto]{\(\scriptstyle \id_{D_k}\)} (ijk-3-2);
\end{tikzpicture},
\end{equation*}
arising from the fact that both cones are initial (that is, both objects are colimits of the above diagram).
\item[(c)] Left and right unitality:\\
For \(i>j\geq 0\), the diagram
\begin{equation*}
\begin{tikzpicture}[bij/.style={above,sloped,inner sep=0.5pt}]
\matrix (idc) [matrix of math nodes, row sep=3em,
column sep=3em, text height=1.5ex, text depth=0.25ex]
{
                  & D_i & \\
D_i \coprod_{D_j} D_j & D_i \coprod_{D_j} D_i & D_j \coprod_{D_j} D_i \\};
\path[->, font=\scriptsize]
(idc-1-2) edge node[auto]{\(\scriptstyle \nabla^i_j\)} (idc-2-2)
(idc-1-2) edge node[bij]{\(\scriptstyle \sim\)} (idc-2-1)
(idc-1-2) edge node[bij]{\(\scriptstyle \sim\)} (idc-2-3)
(idc-2-2) edge node[auto]{\(\scriptstyle \id_{D_i}\coprod_{D_j}\kappa^j_i\)} (idc-2-1)
(idc-2-2) edge node[auto,swap]{\(\scriptstyle \kappa^j_i\coprod_{D_j} \id_{D_i}\)} (idc-2-3);
\end{tikzpicture},
\end{equation*}
commutes.
\item[(d)] Functoriality of units:\\
For \(i>j\geq 0\), the diagram
\begin{equation*}
\begin{tikzpicture}
\matrix (funct) [matrix of math nodes, row sep=3em,
column sep=3em, text height=1.5ex, text depth=0.25ex]
{ D_{i+1} & D_{i+1}\coprod_{D_j} d_{i+1} \\
   D_i & D_{i+1}\coprod_{D_j}\\};
\path[->, font=\scriptsize]
(funct-1-1) edge node[auto]{\(\scriptstyle \nabla^{i+1}_j\)}(funct-1-2)
        edge node[auto,swap]{\(\scriptstyle \kappa_i\)} (funct-2-1)
(funct-2-1) edge node[auto,swap]{\(\scriptstyle \nabla^i_j\)}(funct-2-2)
(funct-1-2) edge node[auto]{\(\scriptstyle \kappa_i\coprod_{D_j} \kappa_i\)} (funct-2-2);
\end{tikzpicture}
\end{equation*}
\end{enumerate}
\end{defn}

\begin{defn}
A morphism  of (pre)categorical \(\G\)-extensions is defined to be a morphism of the underlying \(\G\)-extensions preserving the cocategorical operations \(\nabla^i_j\) and \(\kappa_i\).  
\end{defn}

\begin{defn} Given a category \(C\), a \dfn{strict \(\omega\)-(pre)category internal to \(C\)} is defined to be a globular set \(D^\op:\G^{op}\to C\) together with two specified families of operations \((\nabla^i_j:D_i\to D_i\coprod_{D_j} D_i){i>j\geq 0}\) and \((kappa_i:D_i\to D_{i-1})_{i>0}\) on \(C^\op\) such that the triple \((D,(\nabla^i_j)_{i>j\geq 0},(\kappa_i)_{i>0})\) gives \(C^{op}\) the structure of a (pre)categorical extension.  

By abuse of notation, we will, when the meaning is clear, simply refer to such a triple by its underlying globular set. A morphism \(X\to Y\) of strict \(\omega\)-(pre)categories internal to \(C\) is defined to be a natural transformation \(X\to Y\) such that the induced map \(Y^\op \to X^\op\) respects the operations \(\nabla^i_j\) and \(\kappa_i\) for all \(i>j\geq 0\).   We denote the category of strict \(\omega\)-categories internal to \(C\) by \(\cat{\omega\-Cat}(C)\), or when \(C=\cat{Set}\), simply by \(\cat{\omega-Cat}\).  
\end{defn}

\subsection{$\Theta$ as the initial categorical extension}

We will give a description of \(\Theta\) as a \(\G\)-extension and show that the models for \(\Theta\) are precisely the strict \(\omega\)-categories.  

The first construction, due to Ara in \cite{ara-thesis}, uses his theory of globular envelopes together with a two-step brute-force construction by presentation, which we will recount here:

\begin{enumerate}
\item[(1)] Let \(\Theta_{\operatorname{pcat}}\) denote the globular envelope of the category obtained by formally adjoining the operations \(\nabla^i_j\) and \(\kappa_i\) and taking the quotient by the relations (PC1) and (PC2).
\item[(2)] Let \(\Theta\) denote the globular envelope of the category obtained from \(\Thetap\) by formally identifying the legs of the commutative diagrams in (CC1-CC4).
\end{enumerate}

Then we have the following result following immediately from the universal properties of  and the fact that \(\Theta\) and \(\Thetap\) are, by construction, a categorical extension and a precategorical extension respectively:

\begin{prop}
The canonical functor, natural in \(D\), \[\Mod(\Theta,D)\to \cat{\omega-Cat(D)}\] (respectively, \[\Mod(\Thetap,D)\to \cat{\omega-PCat(D)},\] also natural in \(D\)) is a natural equivalence of categories.
\end{prop}

\subsection{The combinatorial properties of $\Theta$}
We will obtain an explicit definition of the morphisms and objects in \(\Theta\) using the definition from the last section.  First, notice that we have a canonical morphism of \(\G\)-extensions \(\Theta_0\to \Theta\).  We would like to describe this functor using only the axioms for categorical extensions and the globular extension property of \(\Theta_0\).  Then it suffices to describe the hom-sets \(\Hom_\Theta(D_n,S)\) for \(n\geq 0\) and \(S\) any object of \(\Theta\).   We will call a morphism \(f:S\to T\) in \(\Theta\) a \dfn{spinal} monomorphism if it is the image of a morphism in \(\Theta_0\).

Let \(T\) be the object of \(\Theta_0\) (and \(\Theta\)) defined by a globular pattern 
\[D_{i_1}\leftarrow D_{i'_1}\to \dots \leftarrow D_{i'_{k-1}} \to D_{i_k},\] 
and let \(n\geq \heit(T)=\max_{1\leq j\leq k}(i_j)\).  Then there exists a canonical map \(\beta^T_n : \eta_{n(T)} \to T\) defined by the globular pattern \[D_{n}\leftarrow D_{i'_1}\to \dots \leftarrow D_{i'_{k-1}} \to D_{n}\] given by the iterated amalgamation of the appropriate maps \(\kappa^{i_j}_n: D_n \to D_{i_j}\), which we may depict schematically as \[\kappa^{i_1}_n\leftarrow D_{i'_1}\to \dots \leftarrow D_{i'_{k-1}} \to \kappa^{i_k}_n: .\]  Also, suppose we are given a family of nonnegative integers \(i'_j :0\leq j\leq k-1\) and a nonnegative integer \(n\geq i'_j\) for all \(0\leq j\leq k-1\) .  

Then we define the cocomposition operation for an object defined by a globular pattern of the form \[D_{n}\leftarrow D_{i'_1}\to \dots \leftarrow D_{i'_{k-1}} \to D_{n}\] by induction with respect to the special case where \(k=3\).  That is, we define \[\nabla^n_{i'_1,i'_2}:D_n\to D_n \coprod_{D_{i'_1}} D_n \coprod_{D_{i'_2}} D_n\] to be the composite of either leg in the diagram
\begin{equation*}
\begin{tikzpicture}
\matrix (b) [matrix of math nodes, row sep=3em,
column sep=3em, text height=1.5ex, text depth=0.25ex]
{ D_n & D_n \coprod_{D_{i'_1}} D_n \\
D_n\coprod_{D_{i'_2}} D_n &  D_n \coprod_{D_{i'_1}} D_n \coprod_{D_{i'_2}} D_n \\};
\path[->, font=\scriptsize]
(b-1-1) edge node[auto]{\(\scriptstyle \nabla^n_{i'_1}\)} (b-1-2)
        edge node[auto]{\(\scriptstyle \nabla^n_{i'_2}\)} (b-2-1)
(b-2-1) edge node[auto]{\(\scriptstyle \nabla^n_{i'_1}\coprod_{D_{i'_2}} \id_{D_n}\)} (b-2-2)
(b-1-2) edge node[auto]{\(\scriptstyle \id_{D_n} \coprod_{D_{i'_1}}  \nabla^n_{i'_2}\)} (b-2-2);
\end{tikzpicture}.
\end{equation*}    
In fact, the diagram above does commute, which follows from the interchange and unit axioms.  This gives us maps \(\nabla^n_{i'_1,\dots, i'_{k-1}}\) for any choices of \(i'_j\) such that \(i'_j \leq n\) for each \(j\).  Given any object \(T\) of \(\Theta_0\) with height at most \(n\geq 0\), we then define \(c_n^T\) to be the composite map \(D_n\to T\) given by the composite \(\beta^T_n\circ \nabla^n_{i'_1,\dots i'_{k_T}} : D_n\to \eta_n(T) \to T\), and we call it an \dfn{\(n\)-cospine of \(T\)}.  If \(n=\heit(T)\), we call the unique \(n\)-cospine the \dfn{principal cospine of \(T\)}.  We will say that a morphism \(f:S\to T\) in \(\Theta\) is \dfn{cospinal} if \(f\circ c_n^S\) is a cospine for \(T\).

Since the objects of \(\Theta\) represent functors that are \(\omega\)-categories internal to \(\cat{Set}\) (which follows from the fact that \(\Hom_{\Theta}(\cdot,S)\) sends colimits to limits), we see that the principal cospine of an object \(T\) gives the \(T\)-shaped composition map \(\Hom_{\Theta}(T,S)\to \Hom_{\Theta}(D_{\heit(T)},S)\), which gives us the ``total composite cell'' of \(S\) when \(T=S\).  Then the set of \(n\)-globes of the underlying \(\Theta_0\)-model of an object \(S\) in \(\Theta\) is precisely given by the set of pairs consisting of a cospine \(D_n\to T\) together with a spinal monomorphism \(T\to S\), since every globe of \(S\) is uniquely a composite of a subspine \(T\hookrightarrow S\).  That is, \(Hom_{\Theta}(D_n,S)=\coprod_{\heit(T)\leq n}\Hom_{\Theta_0}(T,S)\).  

Since every object \(R\) of \(\Theta\) is uniquely a globular sum of globes \(D_n\), we find that  \[\Hom_{\Theta}(R,S)=\varprojlim_{I_k}\Hom_{\Theta}(D_{n_j},S),\]
which gives us an explicit definition of the hom-sets in \(\Theta\).   

The above discussion implies easily that the following proposition holds:
\begin{prop}\label{spinalfactor} Every morphism \(S\to T\) in \(\Theta\) admits a unique decomposition into a cospinal map followed by a spinal monomorphism.
\end{prop}

We introduce a few definitions that we will use later:

\begin{defn} We call a map \(D_n\to T\) in \(\Theta\) a \dfn{sector} if it is a spinal monomorphism.  When such a map is maximal in the poset \(\overcat{\Theta_0}{T}\), we will call it an \dfn{input sector}.  
\end{defn}
\chapter{A cartesian model structure for weak $\omega$-categories}
\section{Covers, spines, and anodynes}

We will prove some important technical lemmas, which will be absolutely critical to constructing a cartesian-closed model category of weak \(\omega\)-categories.   

\subsection{Rezk Covers}

Ara's theory proves that the objects of \(\Theta\) are in canonical bijection with globular patterns, and further, that the colimit of the globular pattern associated with an object \([t]\in \Theta\), when taken in \(\Theta\) is \([t]\). There is another way to form the colimit of the globular pattern, namely, in the category of presheaves on \(\Theta\).  Then given a globular pattern \(\eta_t\) associated with an object \([t]\in \Theta\), we define its \dfn{spine} to be the colimit of \(h_{\eta_t}\), where \(h_{(\cdot)}\) is the Yoneda embedding.  By the universal property, the spine admits a unique map into the functor \(h_t\) represented by \([t]\).  It is easy to see by induction that the map is injective.  This next definition gives a generalization of the notion of a ``sequential map'' between objects of \(\Delta\).  

\begin{defn}
We say that a map \([s]\to [t]\in \Theta\) is \dfn{spinal} if it sends the spine of \([s]\), written \(\Sp[s]\) into the spine \(\Sp[t]\) of \([t]\).  The monic spinal maps are precisely those maps arising from \(\Theta_0\), so this agrees with our earlier definition.
\end{defn}

\begin{defn} We say that a subpresheaf \(T\subseteq \Theta[t]\) (where \(\Theta[t]\) denotes the functor \(\Hom_{\Theta}(\cdot,[t])\)) is a \dfn{Rezk Cover} provided that:
\begin{enumerate}
\item [(i)] We have \(\Sp[t]\subseteq T\), and
\item [(ii)] The inclusion map \(T\hookrightarrow \Theta[t]\) has the right lifting property with respect to the set of all cospinal maps.    
\end{enumerate}
\end{defn}

\begin{prop}\label{covprops} The following properties hold:
\begin{enumerate}
\item [(i)] Every epimorphism in \(\Theta\) is spinal.
\item [(ii)] The pullback of a cover along a spinal map is a cover.
\item [(iii)] Given two objects \([s]\) and \([t]\), \(S\to \Theta[s]\) and \(T\to \Theta[t]\) two covers, and a pair of spinal maps \([p]\to [s]\) and \([p]\to [t]\), the pullback of the map \(S\times T \hookrightarrow \Theta[s]\times \Theta[t]\) along the map \(\Theta[p]\to \Theta[s]\times\Theta[t]\) is a cover of \([p]\).  
\item [(iv)] The inclusion of the spine is a cover
\item [(v)] The identity map is a cover 
\end{enumerate} 
\end{prop}
\begin{proof}
We leave the proof of these facts to the reader.  
\end{proof}
\subsection{Products of covers are anodyne}
We quickly recall a proposition of Cisinski regarding the behavior of pullbacks under the canonical homotopy colimits with respect to a regular localizer \(W\).  We will state it without proof, and we encourage any uneasy readers to check it in its original source:
\begin{prop}[\cite{cisinski-book}*{Prop. 3.4.46}]\label{fibcolim} Let \((A,W)\) be a small category equipped with a regular \(A\)-localizer. Then given a morphism \(S\to T\) of presheaves on \(A\), we recall the canonical pullback functor
\[\rho_T:A\downarrow T \to \psh{A}, \qquad (a,h_a\to T)\mapsto S\times_T h_a.\]
Then the morphism \(\hocolim^W\rho\to S\), induced by the projection maps \(S\times_T h_a\to S\), belongs to \(W\).
\end{prop}
\begin{lemma}\label{coversweak} Any \(\Theta\)-localizer containing the spine inclusions contains the Rezk covers. 
\end{lemma}
\begin{proof}
Fix a localizer \(W\) of \(\Theta\) containing the spine inclusions, and let \(S\hookrightarrow \Theta[s]\) be a proper \((S\neq \Theta[s])\) cover of an object \([s]\) in \(\Theta\).  By the \(2\)-for-\(3\) property of localizers, it suffices to show that the inclusion \(\Sp[s]\hookrightarrow S\) belongs to \(W\).

Let \(P_S\) denote the category whose objects are the injective spinal maps \(f_p:[p]\hookrightarrow [s]\) that factor through \(S\), and whose morphisms \(([p],f_p)\to ([p'],f_{p'})\) are maps \(g:[p]\to [p']\) such that \(f_p=f_{p'}\circ g\).  We see that this category is isomorphic to the full subcategory of \(\overcat{\Theta}{S}\) spanned by the monomorphisms \(f_p:\Theta_{p}\hookrightarrow S\) such that the composite \(\Theta_{p}\hookrightarrow S\hookrightarrow \Theta[s]\) is a spinal monomorphism.  To every object of \(P_s\), we assign a cartesian rectangle:
\begin{equation*}
\begin{tikzpicture}
\matrix (b) [matrix of math nodes, row sep=3em,
column sep=3em, text height=1.5ex, text depth=0.25ex]
{ Pb_2 & \Sp[s] \\
   Pb_1 &  S \\
	 \Theta[p] & \Theta[s] \\};
\path[->, font=\scriptsize]
(b-1-1) edge (b-1-2)
        edge (b-2-1)
(b-2-1) edge (b-2-2)
        edge (b-3-1)
(b-3-1) edge node[auto]{\(\scriptstyle f_p\)} (b-3-2)
(b-1-2) edge (b-2-2)
(b-2-2) edge (b-3-2);
\end{tikzpicture}
\end{equation*}
However, since \([p]\to [s]\) factors through \(S\), we see that \(Pb_1=h_p\). We also see that \(Pb_2\) is precisely \(\Sp[p]\), since \(f_p\) is spinal and injective. Therefore, we reduce the rectangle above to a cartesian square \(X_{[p]}\):
\begin{equation*}
\begin{tikzpicture}
\matrix (b) [matrix of math nodes, row sep=3em,
column sep=3em, text height=1.5ex, text depth=0.25ex]
{ \Sp[p] & \Sp[s] \\
   \Theta[p] &  S \\};
\path[->, font=\scriptsize]
(b-1-1) edge (b-1-2)
        edge (b-2-1)
(b-2-1) edge (b-2-2)
(b-1-2) edge (b-2-2);
\end{tikzpicture}
\end{equation*}
This square is clearly functorial in \(P_S\).  Since \(\Theta\) is regular squelettique, the localizer \(W\) is necessarily regular, but since \(W\) is regular, \eqref{fibcolim} tells us that for the canonical functor
\[\rho_S:\overcat{\Theta}{S}\to \cellset, \qquad (\theta,h_\theta\to S)\mapsto \Sp[s] \times_S h_\theta,\]
we have that the canonical map \(\hocolim^W \rho_S \to \Sp[s]\) belongs to \(W\).  However, we know that that the inclusion  \(\Sd(S)\hookrightarrow \overcat{\Theta}{S}\) of the full subcategory spanned by the monomorphisms is homotopy cofinal, since \(\Theta\) is skelettique regular, which implies that the natural map \(\hocolim^W \eval[1]{\rho_S}_{\Sd(S)}\to \hocolim^W \rho_S\) belongs to \(W\).  Then we have reduced the problem of a weak equivalence \(\hocolim^W\eval[1]{\rho_S}_{P_S} \to \Sp[s]\) to showing that \(\hocolim^W\eval[1]{\rho_S}_{P_S}\to \hocolim^W \eval[1]{\rho_S}_{\Sd(S)}\) belongs to W.  This fact will certainly follow if \(P_S\) is indeed homotopy-cofinal in \(\Sd(S)\).  We will digress for a few moments:

Since for each object \([p]\to [s]\) of \(P_S\), the lefthand map is a spine inclusion and therefore a weak equivalence, we know by the universal property of homotopy colimits that the canonical map \(\hocolim^W\eval[1]{\rho_S}_{P_S}\to \hocolim^W\eval[1]{\pi_S}_{P_S}\) belongs to \(W,\) where \(\pi_S\) is the obvious forgetful functor \(\overcat{\Theta}{S}\to \cellset\).

Then we would also like to show that the natural map \(\hocolim^W\eval[1]{\pi_S}_{P_S}\to S\) belongs to \(W\).  Similar to the top part of the diagram, we may first reduce this by regularity to the statement that \(\hocolim^W\eval[1]{pi_S}_{P_S}\to \hocolim^W{\pi_S}\) belongs to \(W\), and because \(\Theta\) is skelettique regular, we can further reduce the problem using the cofinality of \(\Sd(S)\) in \(\overcat{\Theta}{S}\), which implies that it suffices to prove that the natural map \(\hocolim^W\eval[1]{\pi_S}_{P_S}\to \hocolim^W\) belongs to \(W\).  As with the top morphism in the diagram, a proof that \(P_S\) is homotopy cofinal in \(\Sd(S)\) will imply that the map in question belongs to \(W\).  

For any injective \(\alpha_q:[q]\hookrightarrow S\), it follows from \eqref{spinalfactor} that there exists a unique factorization of \([q]\hookrightarrow S\hookrightarrow [s]\) into a cospinal map followed by a spinal map \([q]\hookrightarrow [t]\hookrightarrow [s]\), which admits a unique monomorphic lifting \([t]\to S\) making the whole diagram commute.  However, the map \([t]\to S\) belongs to \(P_S\) by inspection, and it is initial in \(\overcat{([q],\alpha_q)}{P_S}\).  This implies that for all \(\alpha_q:[q]\hookrightarrow S\) in \(\operatorname{Sd}_\Theta(S)\), \(\overcat{([q],\alpha_q)}{P_S}\) has a contractible nerve, and therefore \(P_S\hookrightarrow \operatorname{Sd}_\Theta(S)\) is homotopy cofinal.  
\end{proof}
The next lemma establishes our ability to write \([t]\) as an object \([n_t](t_1,\dots,t_n),\) of \(\Delta \wr \Theta\), where each \(t_i\) has height strictly smaller than \([t]\).   
\begin{prop} There exists an isomorphism of categories \(\xi:\Delta \wr \Theta\cong \Theta\).
\end{prop}  
\begin{proof}
Recall that in \eqref{segfun} we defined a functor \(F_\Delta:\Delta\to \Gamma\), by applying the general construction from \eqref{infwreath} to \(\Delta\) equipped with this functor, we constructed the functor \(T_{\Delta,[0],F_\Delta}:\mathbf{N}\to \Cat\)  given by the diagram
\begin{equation*}
\begin{tikzpicture}
\matrix (a) [matrix of math nodes, row sep=3em,
column sep=3em, text height=1.5ex, text depth=0.25ex]
{ \Delta^{\wr 0} & \Delta^{\wr 1} & \Delta^{\wr 2} & \dots & \Delta^{\wr n} & \Delta^{\wr n+1} & \dots \\};
\path[right hook->, font=\scriptsize]
(a-1-1) edge node[auto]{\(\scriptstyle \iota_0\)} (a-1-2)
(a-1-2) edge node[auto]{\(\scriptstyle \iota_1\)} (a-1-3)
(a-1-3) edge node[auto]{\(\scriptstyle \iota_2\)} (a-1-4)
(a-1-4) edge node[auto]{\(\scriptstyle \iota_{n-1}\)} (a-1-5)
(a-1-5) edge node[auto]{\(\scriptstyle \iota_n\)} (a-1-6)
(a-1-6) edge node[auto]{\(\scriptstyle \iota_{n+1}\)} (a-1-7);
\end{tikzpicture},
\end{equation*}
where \(\iota_0:\Delta^{\wr 0}=\ast\to \Delta\) is the functor \([0]:\ast\to \Delta\) classifying the object \([0]\) and \(\iota_{n+1}:\Delta^{\wr n+1}\to \Delta^{\wr n+2}\) is the functor \(\id_\Delta \wr \iota_n\).  Taking the colimit of this diagram, we obtain the category \(\Theta=C(\Delta,[0],F_\Delta)\).  Let \(s:\mathbb{N}\to \mathbb{N}\) be the functor sending \(n\mapsto n+1\).  This functor is clearly cofinal, so the diagrams \(T_{\Delta,[0],F_\Delta}\) and \(s^\ast T_{\Delta,[0],F_\Delta}\) necessarily have isomorphic colimits.  However, \(s^\ast T_{\Delta,[0],F_\Delta}=\Delta \wr T_{\Delta,[0],F_\Delta}\) by construction.  

Then it suffices to show that the functor \(\mathcal{C}\mapsto \Delta\wr \mathcal{C}\) preserves linear colimits.   To see this, note that given a functor \(\Lambda:\mathbf{N}\to \Cat\), we may describe its colimit as the category specified as follows: The set of objects is the quotient of the set of pairs \((x,n)\), where \(x\) is an object of \(\Lambda(n)\) by the equivalence relation \((x,n)\sim(y,m)\) if and only if there exists a natural integer \(p\geq \operatorname{max}(n,m)\) such that \(\Lambda_n^p(x)=\Lambda_m^p(y)\) (where \(\Lambda_i^j:\Lambda(i)\to \Lambda(j)\) is the image of the unique map \(i<j\) in \(\mathbf{N}\)).  We denote the equivalence class of \((x,n)\) by \(\langle x,n\rangle\).  The set of morphisms \(f:\langle x,n\rangle\to \langle y,m\rangle\) is given by the quotient of the set \[\coprod_{i\in \mathbf{N}}\coprod_{(a,i)\in \langle x,n\rangle\\ (b,i)\in \langle y,m\rangle}\Hom_{\Lambda(i)}(a,b)\] modulo the equivalence relation \((f,i)\sim (g,j)\) exactly when there exists \(k\geq\operatorname{max}(i,j)\) such that \(\Lambda_i^k(f)=\Lambda_j^k\).  

Then an object of \(\colim(\Delta\wr T_{\Delta,[0],F_\Delta})\) is an equivalence class of pairs \[\langle [n](x_1,\dots,x_n),i\rangle \] where each \(x_i\) belongs to \(F(i)\), while an object of \(\Delta\wr\colim(T_{\Delta,[0],F_\Delta})\) is of the form \[[n](\langle x_1,i_1 \rangle,\dots,\langle x_n,i_n \rangle).\] 

We see that for any two equivalent families that \[[n]((x_1,i_1),\dots,(x_n,i_n))\sim [n]((y_1,j_1),\dots,(y_n,j_n)),\] since there exists an element \([n]((z_1,k_1),\dots,(z_n,k_n)\) such that \(k_\ell\geq \operatorname{max}(i_\ell,j_\ell)\) for each \(1\leq \ell\leq n\) and such that \(\Lambda_{i_\ell}^{k_\ell}(x_\ell)=\lambda_{j_\ell}^{k_\ell}(y_\ell)\).  

Then by letting \(k=\operatorname{max}_{1\leq \ell\leq n}(k_\ell)\), we see that \[[n]((\Lambda_{k_1}^k(z_1)),\dots,(\Lambda_{k_n}^k(z_k)))\] is also a representative.  

So the map sending the set of pairs \(([n](x_1,\dots,x_n),i)\) to  the set of objects of the form \([n]((y_1,j_1),\dots,(y_n,j_n))\) by the rule \[([n](x_1,\dots,x_n),i)\mapsto [n]((x_1,i),\dots,(x_n,i))\] is compatible with the equivalence relation and also descends to a bijection on equivalence classes.  We leave it to the reader to show that the induced map on \(\Hom\)-sets is also bijective, since the proof is basically identical but notation-heavy.
\end{proof}

\begin{thm}\label{mainthm}
Given a \(\Theta\)-localizer \(W\) containing the spine inclusions, two objects \([s]\) and \([t]\) of \(\Theta\), and two covers \(S\to \Theta[s]\) and \(T\to \Theta[t]\), the map \(S\times T\to \Theta[s]\times \Theta[t]\) belongs to \(W\).  
\end{thm}
\begin{proof}
We define the category \(R_{s,t}\) to be the full subcategory of \(\Sd_{\Theta}(\Theta[s]\times \Theta[t])\) 
spanned by those maps \(\iota_p:\Theta[p]\hookrightarrow\Theta[s]\times \Theta[t]\) such that the composites 
\([p]\to [s]\) and \([p]\to [t]\) are both epimorphic.  For each such \(\iota_p\), we functorially assign a cartesian square
\begin{equation*}
\begin{tikzpicture}
\matrix (b) [matrix of math nodes, row sep=3em,
column sep=3em, text height=1.5ex, text depth=0.25ex]
{ c[p] & S\times T \\
   \Theta[p] &  \Theta[s]\times \Theta[t] \\};
\path[->, font=\scriptsize]
(b-1-1) edge (b-1-2)
        edge (b-2-1)
(b-2-1) edge (b-2-2)
(b-1-2) edge (b-2-2);
\end{tikzpicture}.
\end{equation*}

It follows from \eqref{covprops} that \(c[p]\to \Theta[p]\) is a cover of \([p]\) and therefore a \(W\)-equivalence by the previous lemma.  Then \(\hocolim^W_{R_{s,t}} c[p] \to \hocolim^W_{R_{s,t}} \Theta[p]\) is a \(W\)-equivalence, so by the fact about pullbacks and regular localizers mentioned in the proof of the previous lemma, it suffices to show that \(\hocolim^W_{R_{s,t}}\Theta[p]\to \Theta[s]\times\Theta[t]\) is a \(W\)-equivalence.  To prove this, it suffices to show that \(R_{s,t}\) is homotopy cofinal in \(\Sd_\Theta(\Theta[s]\times \Theta[t])\).  Given a monomorphism \(\alpha_q:\Theta[q]\to \Theta[s]\times \Theta[t]\), we let \(R_{s,t,\alpha_q}=\overcat{([q],\alpha_q)}{R_{s,t}}\).  

From the first description of the category \(\Theta\), we may write \[[x]=[n_x]([x_1],\dots, [x_{n_x}])\] where the height of the \([x_i]\) for \(1\leq i\leq n_x\) is strictly less than the height of \([x]\).  Then for an object \([x]\) in \(\Theta\), we will write \([n_x]\) for the corresponding object of \(\Delta\), and call this the \dfn{\(\Delta\)-collapse} of \([x]\).  Conversely, given an object \([m]\) of \(\Delta\), we write \([m]_0\) for the corresponding object of \(\Theta\).  Both of these associations are functorial, and the first is left adjoint to the second.


Then we write \(Q_{s,t}\) to be the full subcategory of \(Sd_\Delta(\Delta[n_s]\times \Delta[n_t])\) spanned by those maps \(\gamma:\Delta[e]\to \Delta[n_s]\times \Delta[n_t]\) such that the composites \([e]\to [n_s]\) and \([e]\to [n_t]\) are both epimorphic.  Similarly, given \[\alpha_q:\Theta[q]\to \Theta[s]\times \Theta[t],\] we write \[n_{\alpha_q}:\Delta[n_q]\to \Delta[n_s]\times \Delta[n_t]\] for the induced map, and we denote the coslice \[\overcat{([n_q],n_{\alpha_q})}{Q_{s,t}}\] by \[Q_{s,t,\alpha_q}.\]  

Consider the set \(K=\{S,E,SE\}\), and define a graded set of \dfn{southeasterly paths}, \(\mathfrak{Q}=\coprod_i=1^\infty K^i\) with the obvious grading map \(\ell:\mathfrak{Q}\to \mathbf{N}\).  We define a relation on \(\mathfrak{Q}\) where, given a pair of elements \(a,b\in \mathfrak{Q}\), we say that \(a\prec b\) if \(\ell(a)=\ell(b)-1\) and if there exists a natural number \(1\leq i \leq \ell(a)\) such that:
\begin{enumerate}
\item[(i)] The element \(a_i = SE\)
\item[(ii)] We have that \(b_i \neq b_{i+1}\)
\item[(iii)] The elements \(b_i,b_{i+1}\) are in the subset \(\{S,E\}\subset K\).
\item[(iv)] We have that \(a_j=b_j\) for all \(j > i+1\) or \(j<i\).  
\end{enumerate}
Taking the transitive, reflexive closure of this relation gives us a partial order structure on \(\mathfrak{Q}\), since it is clearly antisymmetric.
We define a pair of functions \(d_S,d_T:\mathfrak{Q}\to \mathbf{N}\), where \(d_S\) (respectively \(d_E\)) counts the number of occurrences of the letter \(S\) (respectively the letter \(E\)), including occurrences in \(SE\). For any element \(a\) of \(\mathfrak{Q}\), we call the pair \((d_S(a),d_E(a))\) the terminus of \(a\).  

We can see that \(Q_{s,t}\)  is isomorphic as a poset to the full subposet of \(\mathfrak{Q}\) consisting of the southeasterly paths with terminus \((n_s, n_t)\).  

Let \(\beta_p:\Theta[p] \to \Theta[s] \times \Theta[t]\) be an element of \(R_{s,t}\) lying, naturally, over the map \(n_{\beta_p}:\Delta[n_p]\to \Delta[n_s]\times \Delta[n_t]\) in \(Q_{s,t}\).  By our identification of \(Q_{s,t}\) with the poset of southeasterly paths \(\mathfrak{Q}_{n_s,n_t}\) having terminus \((n_s,n_t)\), we obtain a factorization of any map in \(Q_{s,t}\) into a unique sequence of primitive maps corresponding to the \(\prec\) relation defined above.  Consider the case of a morphism \(f:n_{\beta_p} \to \xi\) of \(Q_{s,t}\) such that under the isomorphism with paths, this witnesses one of the generating relations \(\prec\) such that the path \(n_{\beta_p}\) is obtained by composing a corner of \(\xi\) to a diagonal, that is to say, that \(f\) embeds \(\Delta[n_p]\) as an inner facet of \(\Delta[n_p+1]\) (of course such that \(\xi \circ f = n_{\beta_p}\)).

Then we contend that there is a unique object \(\beta_{p^\prime}:\Theta[p^\prime] \to \Theta[s]\times \Theta[y] \) living over \(\xi\) together with a unique morphism \(\phi:\beta_p \to \beta_{p^\prime}\) lying over \(f\).  

As a map of simplices, \(f\) is the inclusion of an inner face, that is to say, an inclusion \(\delta^i:[n_p]\to [n_p+1]\).   Assume that \(S\) comes first. Let \(k_S=f_S([i-1]([p_1],\dots,[p_{i-1}]))+1\) and \(k_E=f_E([i-1]([p_1],\dots,[p_{i-1}]))+1\).  Then depending on the direction of the corner in \(\xi\) (that is, if we travel south first or east first) we let \([p^\prime]\) be 
\begin{align*}
[n_p+1]([p_1],\dots,[p_{i-1}],[s_{k_S}],[t_{k_E}],[p_{i+1}],\dots, [p_{n_p}]) \intertext{or} [n_p+1]([p_1],\dots,[p_{i-1}],[t_{k_E}],[s_{k_S}],[p_{i+1}],\dots, [p_{n_p}])
\end{align*}
respectively. We obtain totally determined epimorphic maps \([p^\prime]\to [s]\) and \([p^\prime]\to [t]\) by collapsing the the \([1]([t_{k_E}])\) and mapping the \([1]([s_{k_S}])\) on identically to the part of \([s]\) to which it corresponds (and vice versa).  Moreover, it's clear that the induced map \(\beta_{p^\prime}\) into the product is a monomorphism, and even moreover, we see that \([p]\) embeds via \(\phi\) into the \(i\)th simplicial face of \([p^\prime]\) in a way that respects \(\beta_p\) and lives over \(f\).  

Then we show that \(\phi\) is opcartesian, but this is immediate since by construction, \([p^\prime]\) and its structure map are initial with respect to lying over \(\xi\) and accepting a map originating from \(\beta_p\).  Moreover, the uniqueness of the construction shows that it is preserved under composition.

Therefore we see that the collapse functor induces a Grothendieck opfibration from \(R_{s,t}\) to \(Q_{s,t}\).  In fact, the opfibration is stable under coslicing in a way such that the same statement holds for the induced functor  \(R_{s,t,\alpha_q}\) to \(Q_{s,t,\alpha_q}\).  The proof is a tedious check that we leave as an exercise for the reader. 

It follows from \cite{maltsiniotishomotopy}*{2.1.10} that if \(Q_{s,t,\alpha_q}\) (resp. \(Q_{s,t}\)) is weakly contractible, and if \(R_{s,t,\alpha_q}\to Q_{s,t,\alpha_q}\) (resp. \(R_{s,t}\to Q_{s,t}\)) has weakly contractible fibres, then \(R_{s,t,\alpha}\) (resp. \(R_{s,t}\)) is weakly contractible.  However, by \cite{rezk-theta-n-spaces}*{6.12 and 6.13}, we see that \(Q_{s,t,\alpha_q}\) (resp. \(Q_{s,t}\) is weakly contractible.  

In what follows, we will treat the case for \(R_{s,t,\alpha_q}\).  The proof in the case of \(R_{s,t}\) is similar but easier, since it is merely considering the case without the additional constraint of a map \(\alpha_q\).  

Then we would like to prove that the fibres are indeed weakly contractible. Let \(\alpha_{qs}=\pi_s\alpha_q\) and \(\alpha_{qt}=\pi_t\alpha_q\) be the induced maps \(\Theta[q]\to \Theta[s]\) and \(\Theta[q]\to \Theta[t]\).   Then given an object an object \[\Delta[n_q]\xrightarrow{\eta_0} \Delta[e] \xrightarrow{\lambda_{n_s}\times \lambda_{n_t}} \Delta[n_t]\times \Delta[n_s]\] in  \(Q_{s,t,\alpha_q}\), we see that an object of the fibre is given by an object \[\Theta[q]\xrightarrow{\eta} \Theta[p] \xrightarrow{\lambda_s\times\lambda_t} \Theta[s]\times\Theta[t]\] factoring \(\alpha_q\) lying over the point in \(Q_{s,t,\alpha_q}\). 

First, we note that \(\alpha_q\) is given by a pair of families of morphisms \(h^\alpha_{ij}[q_i]\to [s_j]\) for pairs \(i,j\) with \(n_{\alpha_{qs}}(i-1)<j<n_{\alpha_{qs}}(i)\) and \(k^\alpha_{il}[q_i]\to [t_l]\) for pairs \(i,l\) with \(n_{\alpha_{qt}}(i-1)<l<n_{\alpha_{qt}}(i)\) such that the induced maps \[h^\alpha_{ij}\times k^\alpha_{il}: \Theta[q_i]\to \Theta[s_j]\times \Theta[t_l]\] are monic for all appropriate \(i,j,l\).    

Then an object \[\Theta[q]\overset{\eta}{\hookrightarrow} \Theta[p] \xrightarrow{\lambda_s\times\lambda_t} \Theta[s]\times\Theta[t]\] in the fiber over \[\Delta[n_q]\overset{\eta_0}{\to} \Delta[e] \xrightarrow{\lambda_{n_s}\times \lambda_{n_t}} \Delta[n_t]\times \Delta[n_s]\] is given by the data:
\begin{enumerate}
\item[(i)] A family of objects \(([p_1],\dots,[p_e])\) of \(\Theta\) 
\item[(ii)] A family of monomorphisms \(\varepsilon_{ii^\prime}: [q_i]\to [p_{i^\prime}]\) for each pair \(i,i^\prime\) such that \(\eta_0(i-1)< i^\prime\leq \eta_0(i)\)
\item[(iii)] A family of epimorphisms \(f_{i^\prime j}: [p_{i^\prime}]\to [s_j]\) for each pair \(i^\prime,j\) such that \(\lambda_{n_s}(i^\prime-1)<j\leq \lambda_{n_s}(i^\prime)\) (resp. a family of epimorphisms \(g_{i^\prime l}:[p_{i^\prime}]\to [t_l]\) for each pair \(i^\prime,l\) such that \(\lambda_{n_t}(i^\prime-1)<l\leq \lambda_{n_t}(i^\prime)\)) .
\end{enumerate}
satisfying the conditions:
\begin{enumerate}
\item[(a)] The product maps \(f_{ij}\times g_{il}:\Theta[p_i]\to \Theta[s_j]\times \Theta[t_l]\) are injective
\item[(b)] The triple \((\varepsilon_{ii^\prime}, f_{i^\prime j}, g_{i^\prime l})\) gives a factorization \((f_{i^\prime j}\times g_{i^\prime l}) \circ \varepsilon_{ii^\prime}=h^\alpha_{ij}\times k^\alpha_{il}\) for \(\eta_0(i-1)< i^\prime \leq \eta_0(i)\).  
\end{enumerate}

When \(i^\prime\leq \eta_0(0)\) or \(\eta_0(n_q)<j\), the pair \((f_{i^\prime j},g_{i^\prime l})\) specifies 
a unique object of \(R_{s_j,t_l}\).  When \(\eta_0(i-1)< i^\prime \leq \eta_0(i)\), the triple 
\((\varepsilon_{ii^\prime}, f_{i^\prime j}, g_{i^\prime l})\) specifies a unique object of 
\(R_{s_j,t_l,h^\alpha_{ij}\times k^\alpha_{il}}\).  Then we may identify the fibre with 
a product of categories of the form \(R_{s_j,t_l}\) and 
\(R_{s_j,t_l,h^\alpha_{ij}\times k^\alpha_{il}}\), which we call the product decomposition of the fibre.  
In the case of \(R_{s,t}\), the product decomposition of the fibre is simply expressed as a product of categories 
of the form \(R_{s_l,t_k}\).  

We perform well-founded induction on the poset of pairs of natural numbers by letting 
\(A\subseteq \mathbf{N}\times \mathbf{N}\) be the subset of pairs \(a,b\), such that 
for all pairs \([s_0]\) and \([t_0]\) where \(\heit([s_0])=a\) and \(\heit([t_0]=b\), 
the category \(R_{s_0,t_0}\) has a weakly contractible nerve, and for all injective 
maps \(\alpha:\Theta[q_0]\hookrightarrow \Theta[s_0]\times \Theta_[t_0]\), we have 
that \(R_{s_0,t_0,\alpha}\) has a weakly contractible nerve.  Since \(\mathbf{N}\times \mathbf{N}\) 
is well-founded, let \(B=\mathbf{N}\times\mathbf{N} - A\), which by wellfoundedness 
has a minimal element \(a,b\).   Let \([s],[t]\) be a pair such that 
\((\heit([s_]),\heit([t])\) is a minimal element of \(B\) and for which 
the inductive hypothesis fails.  Then the fibres of \(R_{s,t,\alpha}\) (resp. \(R_{s,t}\)) over 
\(Q_{s,t,\alpha}\) (resp. \(Q_{s,t}\)) for some monomorphism \(\alpha:\Theta[q]\to \Theta[s]\times \Theta[t]\) 
admit the aforementioned product decomposition, and since 
\((\heit([s_i]),\heit([t_j]))<(\heit([s]),\heit([t]))\) for any \(i,j\) we see that the 
fibres are products of categories with contractible nerves, and are therefore themselves 
contractible by the continuity of the nerve functor.   Then this implies that \(R_{s,t},\alpha_q\) is weakly contractible for all objects \([s], [t]\) of \(\Theta\) and all appropriate maps \(\alpha_q\).  Then this implies that \(B=\emptyset\), which proves the claim.
\end{proof}

This establishes Theorem \eqref{mainthm}.  We immediately deduce \emph{a fortiori} (since the identity map is a cover) the following corollary:

\begin{cor} If \(W\) is a \(\Theta\)-localizer containing the set of spine inclusions, then for any spine inclusion \(f:\Sp[t]\hookrightarrow \Theta[t]\) and any object \(s\) of \(\Theta\),  the map \(\Theta[s]\times f: \Theta[s]\times \Sp[t]\hookrightarrow \Theta[s]\times \Theta[t]\) belongs to \(W\).
\end{cor}

This corollary may be sharpened using the fact that \(\Theta\) is regular squelettique.

\begin{prop}\label{cartesianness} If \(W\) is a \(\Theta\)-localizer containing the set of spine inclusions, then for any spine inclusion \(f:\Sp[t]\hookrightarrow \Theta[t]\) and any presheaf \(X\) on \(\Theta\), then the map \(X\times f: X\times \Sp[t]\hookrightarrow X\times \Theta[t]\).
\end{prop}
\begin{proof} By \cite{cisinski-book}*{Proposition 8.2.8}, which states that any class of presheaves on a regular squelettique category saturated by monomorphisms and containing the representable presheaves is necessarily the class of all presheaves, it suffices to show that the collection \(C\) of presheaves \(X\) such that \(X\times f\) belongs to \(W\) contains the representable functors and is saturated by monomorphisms.  However, the previous corollary implies the first claim, so it suffices to prove the second. 

To prove that \(C\) is saturated by monomorphisms, notice that any pushout square in which one leg is a monomorphism is a homotopy pushout for the minimal localizer and therefore induces a weak equivalence between the pushouts.  Similarly, any transfinite composition of monomorphisms is a homotopy colimit for the minimal localizer and therefore preserves weak equivalences.  Lastly, closure under retracts follows immediately from the fact that \(W\) is closed under retracts.  This establishes that every presheaf belongs to the aforementioned class \(C\), which establishes the proposition.
\end{proof}

\section{The spine-generated model structure}
Joyal and Cisinski conjectured that a particular a model structure (see \cite{joyal-quategory}) on \(\cellset\), the category of cellular sets, is a model for the weak \(\omega\)-category of weak \(\omega\)-categories, analogous to the Joyal model structure on the category of simplicial sets.  In the absence of a complete description of the ``\(n\)-dimensional inner horns'', they were able to construct the model structure using the technology of localizers developed by Cisinski in \cite{cisinski-book}. Their definition is as follows:

Let \(\mathsf{W}_{\Sp}=\mathsf{W}(S)\) be the \(\Theta\)-localizer generated by the set \(S\) of spine inclusions \(\Sp[t]\hookrightarrow \Theta[t]\) for all \([t]\) in \(\Theta\).   It follows from \eqref{cartesianness} that \(\mathsf{W}(S)\) contains \[\operatorname{cart}(S)=\{X\times f: f\in S\wedge X\in \operatorname{Ob}(\widehat{\Theta})\},\] and therefore, by \cite{cisinski-book}*{Corollary 1.4.19}, we have that \(\mathsf{W}_{\Sp}\) is a \dfn{cartesian \(\Theta\)-localizer}, that is to say, for any morphism \(f:X\to Y\) belonging to \(\mathsf{W}_{\Sp}\) and any presheaf \(Z\) on \(\Theta\), the induced map \(Z\times f\) belongs to \(\mathsf{W}_{\Sp}\).   

By \cite{cisinski-book}*{Theorem 1.4.3}, we see that there exists a unique model structure on \(\widehat{\Theta}\) where the weak equivalences are precisely the elements of \(W_{\Sp}\), and the cofibrations are precisely the monomorphisms.  Since \(\mathsf{W}_{\Sp}\) is cartesian, it follows easily that the aforementioned model structure is cartesian-closed (in the sense that the cartesian product is a left-Quillen bifunctor).  

It follows from the cartesianness of \(\mathsf{W}_{\Sp}\) that \(\mathsf{W}_{\Sp}\)-fibrant objects are are precisely those cellular sets \(X\), which are fibrant in the minimal Cisinski model structure, such that for every spine inclusion \(\Sp[t]\hookrightarrow \Theta[t]\), the induced map \(X^{\Theta[t]}\to X^{\Sp[t]}\) is a trivial fibration.

\subsection{A disproof of the Cisinski-Joyal conjecture}
This isn't the end of the story.  Cisnski and Joyal conjectured in \cite{joyal-quategory} that the fibrant objects in this category model a higher category of weak \(\omega\)-categories.  In fact, this is not so. We sketch below the following explicit counterexample:

\begin{thm} Let \([1]\to [1](G_2)\) be the map of strict \(\omega\)-categories obtained from the inclusion \(\ast \to G_2\).  This map is a strict \(\omega\)-equivalence of strict \(\omega\)-categories.  However, the image of this map under the \(\Theta\)-nerve does not belong to \(\mathsf{W}_{\Sp}\).
\end{thm}
\begin{proof}We first note that \(\mathfrak{N}([1](X))\) is necessarily \(\mathsf{W}_{\Sp}\)fibrant for any strict \(\omega\)-category \(X\).  To see this, we first notice that such an object is minimally fibrant, since \(\mathfrak{N}([1](X))^J\cong \mathfrak{N}([1](X)^{G_2})\), which follows from the fact that the category of strict \(\omega\)-categories is cartesian-closed and embeds fully and faithfully in \(\cellset\).  However, we notice that \([1](X)^G_2\) is isomorphic to \([1](X)\), since a functor \(Z\to [1](X)^{G_2}\) is given by precisely the data of a natural \(1\)-isomorphisms between two functors \(Z\to [1](X)\).  However, since the only \(1\)-isomorphisms in this category are identities, there is exactly one such functor for each functor \(Z\to [1](X)\), which means that they are isomorphic.  

It is also easy to see that for any strict \(\omega\)-category \(X\), the map induced by a spine inclusion \(\Sp[t]\hookrightarrow \Theta[t]\), that is, \(\mathfrak{N}(X)^{\Theta[t]}\to \mathfrak{N}(X)^{\Sp[t]}\), is an isomorphism, which again follows from the fact that the category of strict \(\omega\)-categories is cartesian-closed, together with the characterization of the nerves of strict \(\omega\)-categories as those presheaves sending the objects \([t]\), which are globular sums in \(\Theta\), to globular products in the category of sets.  

Since \(\mathfrak{N}([1])\to \mathfrak{N}([1](G_2))\) is a map between fibrant-cofibrant objects, it is necessarily a strong deformation retract, but we know that \((-)\times J\) is a functorial cylinder, since the map \(J\to e\) is a trivial fibration.  If the map is a strong deformation retract, then this can be exhibited by means of a \(J\)-homotopy between the maps in question.

However, such a homotopy \(\mathfrak{N}([1](G_2))\times J \to \mathfrak{N}([1](G_2))\) necessarily lies in the image of the \(\Theta\)-nerve, which means that it corresponds exactly to a natural isomorphism between the functors maps.  However, we know that such a map cannot exist, because \([1](G_2)\) contains no nontrivial isomorphisms and is not isomorphic to \([1]\).  Therefore, the map on nerves cannot be a \(\mathsf{W}_{\Sp}\)-equivalence.
\end{proof}
In particular, this implies that the model structure generated by this localizer cannot be repaired without adding new weak equivalences.  That is to say, there is absolutely no way to fix it by adjusting only the fibrations and cofibrations, which answers Joyal's original conjecture in the negative.

However, all hope is not lost.  The question, then, is how to enlarge \(\mathsf{W}_{\Sp}\). The two obvious ways are to attempt to stabilize the localizer under suspension and find a reasonable set of generators, or to find a new cylinder functor whose homotopy relation encodes a weaker notion of equivalence than strict isomorphism. With youthful na\"ivet\'e, we run headlong into a trap by asking what happens if we simply stabilize \(\mathsf{W}_{\Sp}\) under suspension.

\section{Stable $J$-homotopy and Stably isofibrant cellular sets}
In this section, we will give a stabilized model structure using the theory discussed in \eqref{weakenrichment}.  Given its close relationship with Rezk's theory of \(\Theta\)-enrichment, it will come as no surprise to the reader that the approach we take will be equivalent to the limiting case \((n=\infty)\)of the definition of a weak \(n\)-category discussed by Rezk in \cite{rezk-theta-n-spaces}.  We will also prove that, to end up with a definition of a weak \(\omega\)-category that extends the homotopy theory of strict \(\omega\)-categories discussed in \cite{lmw}, we must localize still further.

\subsection{The $\Theta$-localizer $\W_{\operatorname{StabIso}}$}
We let \(\W_0=\W_\Sp\) be the na\"ive Cisinski-Joyal \(\Theta\)-localizer.  Since \(\Theta \cong \Delta\wr\Theta\), we will show that \(\W_0\) is strongly generated by the set of spine inclusions \(\Sp[t]\hookrightarrow \Theta[t]\) together with the single map \(j:J\to e\).  
\begin{lemma}\label{spinescontaincores} 
The \(\Theta\)-localizer \(\W_0\), viewed as a \(\Delta\wr \Theta\)-localizer, contains the Segal core inclusions. 
\end{lemma} 
\begin{proof}
First, we note that the maps of the form \(\Delta_1[\Sp[s]]\hookrightarrow \Delta_1[\Theta[s]]\) belong to \(\W_0\) by merit of the fact that \(\Delta_1[\Sp[s]]\) is precisely the spine of \(\Delta_1[\Theta[s]]\).

Then we notice that for any two objects \(s,t\) of \(\Theta\), we can show that the inclusion \[\Delta_1[\Sp[s]]\coprod_{\{0\}} \Delta_1[\Sp[t]]\hookrightarrow \Delta_1[\Theta[s]]\coprod_{\{0\}} \Delta_1[\Theta[t]\] belongs to \(\W_0\) by merit of the fact that it is a composite of pushouts of trivial cofibrations.  

By induction, this implies that for any family \(t=(t_1,\dots,t_n)\) of objects in \(\Theta\), if we let \(\Sp[t]=(\Sp[t_1],\dots,\Sp[t_n])\), the canonical map \[\Sc_n[\Sp[t]]\hookrightarrow \Sc_n[t]\] also belongs to \(\W_0\).  However, it is easy to see that \(\Sc_n[\Sp[t]]=\Sp[\Delta_n[t]]\), and therefore that the composite of the two maps  \[\Sc_n[\Sp[t]]\hookrightarrow \Sc_n[t]\hookrightarrow \Delta_n[t]\] is a spine inclusion, which also belongs to \(\W_0\).  Therefore, by \(3\)-for-\(2\), it follows that \(\Sc_n[t]\hookrightarrow \Delta_n[t]\) belongs to \(\W_0\).  
\end{proof}
Then by the second assertion in \eqref{weakcats}, we obtain the following corollary:
\begin{cor}The localizer \(\W_0\) is strongly generated by the small set \[\Sp \cup \{j:J\to e\},\] where \(\Sp\) denotes the set of all spine inclusions \(\Sp[t]\hookrightarrow \Theta[t]\).  
\end{cor} 
Then we define an increasing sequence of localizers using the isomorphisms \[\Theta\cong \Delta\wr \Theta \cong \dots \cong \Delta^{\wr n}\wr \Theta.\]
\begin{defn} For \(n>0\), using the main theorem of \eqref{weakenrichment}, we define \(\W_n\),  the \(n\)-suspended Cisinski-Joyal localizer by the formula \[\W_n=(\W_{n-1})_\wr.\]
\end{defn}
\begin{note}We will denote the n-fold iterate of the suspension \(\Delta_1[-]\) by \((\Delta_1)^n[-]\) (or sometimes also by \(D_n[-]\)), where \((\Delta_1)^0[-]\) denotes the identity functor.  We make note of this, since the notation could also mean the suspension along the \(n\)-fold cartesian power of the simplicial set \(\Delta_1\).  We will denote that functor instead by \((\Delta_1^n)[-]\).
\end{note}
\begin{lemma} For all \(n\geq 0\), \(\W_n\subseteq \W_{n+1}\).  
\end{lemma}
\begin{proof}We see from \eqref{weakenrichment} that \(\W_{n+1}\) is strongly generated by \[S_{n+1}=\Sc \cup \{j\} \cup \Delta_1[S_n],\] so by running the proof of \eqref{spinescontaincores} in reverse, we see that \(\W_n\) contains the full set of spines \(\Sp\).  Then the only strong generators in this class not belonging to \(\Sp\) are the maps of the form \(\Delta_1^k[j]:\Delta_1^k[J]\to \Delta_1^k[e]\) for all \(0\leq k\leq n+1\).  This proves the claim.
\end{proof}
\begin{defn}We define a new \(\Theta\)-localizer \[\W_{\operatorname{StabIso}}=\bigcup_{n\geq 0} \W_n,\] and we call the model structure it generates the \dfn{isostable Joyal model structure} for cellular sets.  
\end{defn}
\begin{prop} The \(\Theta\)-localizer \(\W_{\operatorname{StabIso}}\) is cartesian and strongly generated by the set of maps \[S_\omega=\Sp\cup \{(\Delta_1)^n[j]:n\geq 0\}.\]
\end{prop}
\begin{proof} This follows from the previous lemma together with \eqref{weakenrichment}.
\end{proof}
Then we obtain the following corollary:
\begin{cor} The \(\Theta\)-localizer \(\W_{\operatorname{StabIso}}\) is stable under the weak enrichment process.  That is, \((\W_{\operatorname{StabIso}})_\wr=\W_{\operatorname{StabIso}}\).
\end{cor}
This gives \(\W_{\operatorname{StabIso}}\) the following stability property:
\begin{prop}If \(f:A\to B\) is an arrow belonging to \(\W_{\operatorname{StabIso}},\) then we have that \(\Delta_1[f]:\Delta_1[A]\to \Delta_1[B]\) also belongs to \(\W_{\operatorname{StabIso}}\).  
\end{prop}
\begin{proof}Since the functor \(\Delta_1[-]\) is a parametric left adjoint that preserves cofibrations, the left-adjoint factor \(D_1:\cellset \to \overcat{\Delta_1[\emptyset]}{\cellset}\) is a cofibration-preserving left-adjoint.  It follows from the main theorem of \eqref{weakenrichment} that \(\Delta_1[J\times X] \to \Delta_1[X]\) belongs to \(\W_{\operatorname{StabIso}}\) for every cellular set \(X\).

Since \(\Delta_1[J\times X]\to \Delta_1[X]\) belongs to \(\W_{\operatorname{StabIso}}\) for every cellular set \(X\), the associated maps between bipointed objects \(D_1[J\times X]\to D_1[X]\) are weak equivalences for the coslice model structure on \(\overcat{\Delta_1[\emptyset]}{\cellset}\), which means that \(D_1^{-1}(\overcat{\Delta_1[\emptyset]}{\W_{\operatorname{StabIso}}})\) is a weakly saturated class of maps such that every object admits a cylinder, so by \cite{cisinski-book}*{Proposition 1.4.13}, it is a \(\Theta\)-localizer.  Moreover, it contains \(S_\omega\), which is a set of strong generators for \(\W_{\operatorname{StabIso}}\), so it contains \(\W_{\operatorname{StabIso}}\), which implies the proposition.   
\end{proof}
\subsection{Stable $J$-homotopy and its failure to model $\omega$-equivalence}\label{stablefailure}
We will now proceed to prove that \(\W_{\operatorname{StabIso}}\) still fails to capture the full notion of \(\omega\)-equivalence of strict \(\omega\)-categories. To do this, consider the following: Let \(J_i\) be the unique contractible groupoid containing \(i+1\) objects. This associaction determines an obvious functor \(\Delta\to \cat{Gpd}\), where  \([i]\mapsto J_i\).  Embedding the \(J_i\) into the category of strict \(\omega\)-categories, we define \(E_i=Q(J_i)\), where \(Q\) is the functor sending a strict \(\omega\)-category \(X\) to its universal polygraph resolution.  This gives a cosimplicial object in the category of strict \(\omega\)-categories.
\begin{lemma} The presheaf \(P=\mathfrak{N}_\Theta(E_1)\) is not trivially fibrant in the category of cellular sets. 
\end{lemma}
\begin{proof} Since \(E_1\) is fibrant and contractible in the category of strict \(\omega\)-categories, we can choose two maps \(\eta_0, \eta_1: \partial D_1 \to E_1\) such that \(s(\eta_0) = t(\eta_1)\), and \(t(\eta_0) = s(\eta_1)\), where \(s(\eta_0)\neq s(\eta_1)\).  Using its fibrancy and contractibility, we extend these to maps \(D_1 \to P\), and by merit of the conditions earlier imposed, we obtain a map \(\sigma:\Lambda^1[2]\to P\).  However, this must lift to a map \(\sigma^\prime:\Delta[2] \to P\) because \(P\) is nerve of a strict \(\omega\)-category. Moreover, \(d_1(\sigma^\prime)\) is the composite of two non-identity maps in a strict \(\omega\)-category freely generated by polygraph, and is therefore a nondegenerate edge.  Then consider the extension of \(\sigma\) to \(\partial\Delta[2]\) with a degenerate edge at the face opposite the first vertex.  However, this clearly does not admit a lift to a full simplex, since it would contradict the uniqueness of the original simplex we found.  Then it is not trivially fibrant.  
\end{proof}
If it is fibrant for the stabilized Cisinski-Joyal model structure, then it cannot be weakly contractible, which implies that the model structure is still wrong and needs to be localized still further.
\begin{thm} The cellular set \(P\) is not weakly contractible with respect to \(\W_{\operatorname{StabIso}}\).
\end{thm}
\begin{proof} Because \(\W_{\operatorname{StabIso}}\) is a cartesian localizer generated by a set of monomorphisms \(S\), showing that \(P\) is fibrant amounts to proving that \(P^f\) is a trivial fibration for all \(f \in S\). If \(f\in S\) is a spine inclusion, we see that \(P^f\) is an isomorphism, since \(P\) is the nerve of a strict \(\omega\)-category.  We also know that the map is a trivial fibration when \(f=j:\Theta[0] \hookrightarrow J\).  It suffices, then, to prove the cases where \(f=D_k[j]:D_k\hookrightarrow D_k[J]\).  That is, it suffices to show that \(P\) has the right lifting property with respect to all maps of the form \[D_k[J] \times \partial\Theta[t] \cup D_k \times \Theta[t] \hookrightarrow D_k[J] \times \Theta[t].\]  If we consider what it means to give a map \(D_k[J]\times \Theta[t] \to P\), this involves giving a map from each shuffle of \(D_k[J]\) with \(\Theta[t]\) and asking them to agree on intersections.
Explicitly, if \([t]=[n_t]([t_1],\dots, [t_{n_t}])\), 
the shuffles are exactly the objects 
\(U_i=[n_t+1](\dots, [t_i], D_{k-1}[J],[t_{i+1}],\dots)\).  
However, since \(P\) contains no strictly invertible non-identity cells (as it is obtained from an \(\omega\)-polygraph), it follows that a map \(U_i \to P\) would factor through \(V_i=[n_t+1](\dots, [t_i], D_{k-1},[t_{i+1}],\dots)\).  However, these are precisely the shuffles of 
\(D_{k}\) with \(\Theta[t]\).   Then we can simply choose the lifting in the obvious way, by mapping each \(U_i\) back onto \(V_i\) and back through the map \(D_k\times \Theta[t]\).  This proves the theorem.
\end{proof}
\section{A cartesian presentation for weak $\omega$-categories}

Consider the inclusion \(\iota:\Theta[0]\to P\) of one of the vertices into \(P\).   Recall that the category of simplicial presheaves on \(\Theta\) is canonically simplicially enriched by the formula \(\HMap_\Delta(X,Y)=(X^Y)[0]\), and that this is well-behaved with respect to the model structure.  

\begin{lemma}
If X is stably Segal, and if \(\HMap_\Delta(\iota \times D_n, X): \HMap_\Delta(P\times D_n,X) \to \HMap_\Delta(\Theta[0]\times D_n, X)\) is a trivial fibration for all \(n\geq 0\), then so too is \(X^\iota:X^P\to X\). 
\end{lemma}
\begin{proof}
To begin with, \(X^\iota\) is a trivial fibration if and only if \(\HMap_\Delta(\Theta[s],X^\iota)\) is a trivial fibration for all discrete representables \(\Theta[s]\), which is equivalent to asking that \(\HMap_\Delta(\iota, X^{\Theta[s]})\) is a trivial fibration for all choices of \([s]\).  However, we have trivial fibrations, \(\HMap_\Delta(I,X^{\Theta[s]}) \to \HMap_\Delta(I,X^{\Sp[s]})\) and \(\HMap_\Delta(\Theta[0],X^{\Theta[s]}) \to \HMap_\Delta(\Theta[0],X^{\Sp[s]})\) because \(X\) is stably Segal.  However, \(\Sp[s]\) decomposes as the colimit over a globular sum diagram whose maps are cofibrations between disks, so \(X^{\Sp[s]}\) is isomorphic to the limit over a coglobular product diagram, where the maps are fibrations between objects of the form \(X^{D_i}\).  Then this limit is a homotopy limit of a globular product diagram with entries \(X^{D_i}\).  Then \(\HMap_\Delta(p,X^{\Sp[s]}) \) is a trivial fibration if and only if \(\HMap_\Delta(\iota,X^{D_i})\) for each appropriate \(D_i\).  The result follows by adjunction.
\end{proof}

Using this, we may form the \(\Theta\)-localizer \(\W_1=\W(\Sp \cup \{\iota \times D_n\}_{n\in \mathbf{N}})\).  The set of arrows \(\Sp \cup \{\iota \times D_n\}_{n\in \mathbf{N}} \cup \{j\}\) forms a set of strong regular generators for \(\W_1\).  As before, we form the suspension stabilization of this localizer using the above class of strong generators and call this localizer \(\W_\omega\).  We see that \(\W_\omega\) is cartesian and suspension stable.  

Ultimately, we do not end up with a presentation as nice as the one in \cite{rezk-theta-n-spaces} or the equivalent one in \cite{ara-new} for \(n\)-categories with \(n\) finite.  This stems from our current inability to prove a version of Rezk's theorem showing that \(J\) is a model for a homotopy equivalence.  This stalls us at the step where we reduce to showing that an object local with respect to \(\iota\) is local with respect to \(\iota\times D_n\) for each \(n\geq 0\).  We state a conjecture, which would give a nicer presentation:
\begin{conj}\label{genconjecture} If X is stably Segal, and if \(\HMap_\Delta(\iota, X): \HMap_\Delta(P,X) \to \HMap_\Delta(\Theta[0], X)\) is a trivial fibration, then so too are \(X^\iota:X^P\to X\) and \(X^j:X^J\to X\).  
\end{conj}
The resulting presentation would then only depend on suspending the strong generators \(\Sp \cup \{\iota\}\), which would give a presentation that is \emph{as nice} as the ones described for \((\infty,n)\)-categories when \(n\) is finite.  

However, we will now describe a completely different approach to the construction of a similar model structure, which may be equivalent.  If it is, then it gives us a powerful new way to deal with the theory.



\begin{bibdiv}
\begin{biblist}

\bib{ara-thesis}{thesis}{
	author={Ara, D.},
	title={Sur les \(\infty\)-groupo\"{i}des de Grothendieck},
	organization={Universit\'{e} Paris Did\'{e}rot (Paris 7)},
	date={2010}
	}
\bib{ara-new}{article}{
  author={Ara, D.},
  title={Higher quasi-categories vs higher Rezk spaces},
	eprint={arXiv:1206.4354 [math.AT]}
}

\bib{berger-cellular-nerve}{article}{
  author={Berger, C.},
  title={A cellular nerve for higher categories},
  journal={Adv. Math.},
  volume={169},
  date={2002},
  number={1},
  pages={118--175},
  issn={0001-8708},
  review={\MR {1916373 (2003g:18007)}},
}

\bib{berger-iterated-wreath}{article}{
  author={Berger, C.},
  title={Iterated wreath product of the simplex category and iterated loop spaces},
  journal={Adv. Math.},
  volume={213},
  date={2007},
  number={1},
  pages={230--270},
  issn={0001-8708},
  review={\MR {2331244 (2008f:55010)}},
}

\bib{cisinski-book}{book}{
author={Cisinski, D.-C.},
title={Les pr\'efaisceaux comme mod\`eles des types d'homotopie},
publisher={Soc. Math. France},
date={2006},
series={Ast\'erisque},
volume={308},
}

\bib{cisinski-decalage}{article}{
	author={Cisinski, D.-C.},
	author={Maltsiniotis, G.},
	title={La cat�gorie \(\Theta\) de Joyal est une cat�gorie test},
	journal={J. Pure Appl. Algebra},
	date={2011},
	volume={215},
	pages={962--982},
}

\bib{joyal-theta-note}{article}{
  author={Joyal, A.},
  title={Disks, duality, and Theta-categories},
  date={1997-09},
  status={preprint},
}

\bib{joyal-quategory}{article}{
  author={Joyal, A.},
  title={The theory of quasi-categories and its applications},
  conference={ title={Advanced Course on Simplicial Methods in Higher Categories}, 
	address={Bellaterra, Spain},
	book={ 
	title={Course notes for the Advanced Course on Simplicial Methods in Higher Categories}, 
	address={Bellaterra, Spain}, 
	organization={Centre de Recerca Matem�tica}, 
	},
	note={Available online},},
  date={2008},
  pages={150--496},
	}
	
\bib{htt}{book}{
	author={Lurie, J.},
	title={Higher Topos Theory},
	publisher={Princeton University Press},
	date={2006},
	series={Annals of Mathematics Studies},
	volume={170},
	}
	
\bib{lmw}{article}{
	author={Lafont, Y.},
	author={M\'etayer, F.},
	author={Worytkiewicz, K.},
	title={A folk model structure on omega-cat},
	journal={Adv. Math.},
	date={2007},
	volume={224},
}

\bib{maltsiniotishomotopy}{book}{
author={Maltsiniotis, G.},
title={La th\'eorie de l'homotopie de Grothendieck},
publisher={Soc. Math. France},
date={2005},
series={Ast\'erisque},
volume={301},
}

\bib{rezk-segal-spaces}{article}{
  author={Rezk, C.},
  title={A model for the homotopy theory of homotopy theory},
  journal={Trans. Amer. Math. Soc.},
  volume={353},
  date={2001},
}

\bib{rezk-theta-n-spaces}{article}{
  author={Rezk, C.},
  title={A Cartesian presentation of weak \(n\)-categories},
  journal={Geom. Topol.},
  volume={14},
  date={2010},
  number={1},
  pages={521--571},
  issn={1465-3060},
  review={\MR {2578310}},
  doi={10.2140/gt.2010.14.521},
}

\bib{rezk-theta-n-spaces-correction}{article}{
  author={Rezk, C.},
  title={Correction to "A Cartesian presentation of weak \(n\)-categories"},
  journal={Geom. Topol.},
  volume={14},
  date={2010},
  number={1},
  pages={2301-2304},
  doi={10.2140/gt.2010.14.2301},
}

\bib{street}{article}{
	author={Street, R.},
	title={The petit topos of globular sets},
	journal={J. Pure Appl. Algebra},
	date={2000},
	volume={154},
	pages={299--315},
	}


\end{biblist}
\end{bibdiv}
\end{document}